\documentclass{siamart190516}
\usepackage{amsmath}
\usepackage{float}

\usepackage[font=small,labelfont=bf]{caption}
\usepackage{afterpage}
\usepackage{amsfonts}

\usepackage[labelformat=simple]{subcaption}

\usepackage{amssymb}
\usepackage{graphicx}
\usepackage{amssymb}
\usepackage{amsfonts}
\usepackage{url}
\usepackage{epstopdf}
\usepackage{color}
\usepackage{tikz-cd}
\usepackage{bm}
\usepackage{multirow}
\usepackage{rotating}

\newcommand{\bsub}{\begin{subequations}}
\newcommand{\esub}{\end{subequations}$\!$}

\newcommand{\sumcapone}{\sum_{j=1}^{\scriptscriptstyle N}}
\newcommand{\sumcaptwo}{\sum_{\stackrel{i=1}{i\neq j}}^{\scriptscriptstyle N}}
\newcommand{\sumcapthree}{\sum_{\stackrel{k=1}{k\neq i}}^{\scriptscriptstyle N}}

\newcommand{\eps}{{\varepsilon}}
\newcommand{\X}{{\bf X}}
\newcommand{\K}{{\mathcal K}}

\newcommand{\pa}{{\partial {\mathcal B}}}
\newcommand{\pae}{{\partial_{e} {\mathcal B}}}

\newcommand{\x}{{\bf{x}}}

\newcommand{\vv}{{\bf{v}}}
\newcommand{\vc}{{\bf{C}}}

\newcommand{\vb}{{\bf{b}}}
\newcommand{\y}{{\bf{y}}}
\newcommand{\bxi}{{\boldsymbol{\xi}}}

\newcommand{\R}{{\mathbb{R}}}

\newcommand{\p}{\prime}
\newcommand{\PT}{{\Gamma}}
\newcommand{\DT}{{\bf p}}

\newcommand{\ceff}{C_{\rm eff}}
\newcommand{\keff}{k_{\rm eff}}

\newcommand{\Keff}{{\mathcal K}_{\rm eff}}

\def\KK{{\bf K}}
\def\MM{{\bf M}}

\newtheorem{prop}{Proposition}

\renewcommand\thesection{\arabic{section}}
\renewcommand\thesubsection{\thesection.\arabic{subsection}}
\renewcommand\thesubsubsection{\thesubsection.\arabic{subsubsection}}
\renewcommand\theequation{\thesection.\arabic{equation}}

\graphicspath{{figures/}}

\title{The Effective Reactivity for Capturing Brownian Motion
  by Partially Reactive Patches on a Spherical Surface}

\author{Denis S. Grebenkov \thanks{Laboratoire de Physique de la
    Mati\'{e}re Condens\'{e}e, CNRS -- Ecole Polytechnique,
    Institut Polytechnique de Paris, 91120 Palaiseau, France. email:
    denis.grebenkov@polytechnique.edu}, \and Michael
  J. Ward \thanks{Department of Mathematics, University of British
    Columbia, Vancouver, B.C., Canada, V6T 1Z2. email:
    ward@math.ubc.ca (corresponding author)}}

\date{\today}

\begin{document}

\label{firstpage}
\maketitle

\baselineskip=12pt

\begin{abstract}
  We analyze the trapping of diffusing ligands, modeled as Brownian
  particles, by a sphere that has $N$ partially reactive boundary
  patches, each of small area and arbitrary shape, on an
  otherwise reflecting boundary. For such a structured target, the
  partial reactivity of each boundary patch is characterized by a
  Robin boundary condition, with a local boundary reactivity
  $\kappa_i$ for $i=1,\ldots,N$. For any spatial arrangement of
  well-separated patches on the surface of the sphere, the method of
  matched asymptotic expansions is used to derive explicit results for
  the capacitance $C_{\rm T}$ of the structured target, which is valid
  for any $\kappa_i>0$.  This target capacitance $C_{\rm T}$ is
  defined in terms of a Green's matrix, which depends on the spatial
  configuration of patches, the local reactive capacitance
  $C_i(\kappa_i)$ of each patch and another coefficient that depends
  on the local geometry near a patch. The analytical dependence of
  $C_{i}(\kappa_i)$ on $\kappa_i$ is uncovered via a spectral
  expansion over Steklov eigenfunctions.  For circular patches, the
  latter are readily computed numerically and provide an accurate
  fully explicit sigmoidal approximation for $C_{i}(\kappa_i)$.  In
  the homogenization limit of $N\gg 1$ identical uniformly-spaced
  patches with $\kappa_i=\kappa$, we derive an explicit scaling law
  for the effective capacitance and the effective reactivity of the
  structured target that is valid in the limit of small patch area
  fraction. From a comparison with numerical simulations, we show that
  this scaling law provides a highly accurate approximation over the
  full range $\kappa>0$, even when there is only a moderately large
  number of reactive patches.
\end{abstract}

\section{Introduction}\label{all:intro}

Diffusive search for multiple targets is a key process in physics,
chemistry, and biology
\cite{House,Murrey,Redner,Schuss,Metzler,Holcman15,Lindenberg,Grebenkov,Dagdug}.
For instance, signal transduction between neurons involves several
diffusive search processes
\cite{Alberts,Berridge03,Guerrier18,Reva21}: (i) calcium ions that
diffuse in the extracellular space and search for ionic channels on
the plasma membrane to respond to an electrical signal; (ii) after
entering the synaptic bouton, calcium ions search for calcium-sensing
proteins on the vesicles filled with neurotransmitters to initiate
their release into the inter-neuronal space; (iii) once released, the
neurotransmitters diffusively search for suitable binding sites on the
plasma membrane of a neighboring neuron.  Likewise, small metabolites
and various proteins such as transcription factors or histons search
for nuclear pores on the plasmic membrane of the cell nucleus
\cite{Alberts,Lauffenburger}.  From a chemical perspective,
heterogeneous catalysis often involves porous catalysts with
heterogeneous distribution of active sites on an otherwise passive
(inert) boundary \cite{House,Coppens99,Filoche08,Punia21}.  In both
applications, one generally deals with diffusive search for small
targets or reactive patches that are distributed on an otherwise
reflecting surface.

With this biophysical motivation, we consider the canonical problem of
the trapping of ligands, modeled as diffusing Brownian particles, by
the boundary $\partial{\mathcal B}$ of a 3-D simply-connected and
bounded domain ${\mathcal B}$ that has many small partially reactive
sites.  The boundary is assumed to be smooth, reflecting,
  and partially covered by $N$ small reactive patches
$\partial {\mathcal B}_i$.  The steady-state concentration
${\mathcal U}$ of diffusing ligands satisfies the mixed Neumann-Robin
boundary value problem (BVP) \bsub \label{mfpt:ssp0}
\begin{align}
  \Delta {\mathcal U} & = 0 \,, \quad \X \in \R^3\backslash {\mathcal B} \,,
  \\
   D\partial_{n} {\mathcal U} + \K_i {\mathcal U} & = 0\,, \quad \X \in
     \partial {\mathcal B}_i \,, \quad i=1,\ldots,N \,, \\
  \partial_{n} {\mathcal U} & = 0 \,, \quad \X \in \partial {\mathcal B}_r
                              = \partial{\mathcal B} \backslash 
\overline{(\partial {\mathcal B}_1 \cup \ldots \cup \partial {\mathcal B}_N)}\,,\\
  {\mathcal U} & \sim  {\mathcal U}_{\rm \inf}\left(1 -\frac{{\mathcal C}_{\rm T}}
                 {|\X|} +    {\mathcal O}(|\X|^{-2})\right)
  \,,  \quad \mbox{as} \quad |\X|\to \infty \,, \label{mfpt:ssp0_c}
\end{align}
\esub where
$D>0$ is a constant diffusivity, ${\mathcal U}_{\rm \inf}$ is a
constant concentration imposed at infinity, $\Delta$ is the Laplacian
in the dimensional coordinate $\X$, and $\partial_n$ is the outward
normal derivative to $\partial {\mathcal B}$, directed into
${\mathcal B}$. Each reactive boundary patch $\partial {\mathcal B}_i$
of diameter $2L_i$, with reactivity $\K_i>0$, is assumed to be
simply-connected with a smooth boundary, but with an otherwise
arbitrary shape.  From the divergence theorem, the coefficient
${\mathcal C}_{\rm T}$ in the far-field (\ref{mfpt:ssp0_c}) is related to
${\mathcal U}$ by the identity
\begin{equation}\label{eq:j_origin}
  {\mathcal C}_{\rm T} = -\frac{1}{4\pi {\mathcal U}_{\rm \inf}}
  \int\limits_{\partial {\mathcal B}} \partial_n {\mathcal U} \, ds
= \frac{J}{4\pi D {\mathcal U}_{\rm \inf}} \,,
\end{equation}
where $J$ is the total flux of particles reacted onto patches.  This
coefficient ${\mathcal C}_{\rm T}$ is interpreted as the capacitance of the
heterogeneous partially reactive boundary $\partial {\mathcal B}$.

Our focus will be the spherical domain ${\mathcal B} = \{ \X\in \R^3
\, \vert \, |\X| \leq R\}$, for which we will derive an
analytical formula for the capacitance ${\mathcal C}_{\rm T}$ in the limit
of small patch radius $\eps={L/R}\ll 1$, where $L=\max_{i} \{L_i\}$.
For a homogeneous reactivity $\K$ (i.e., a single patch covering the
whole spherical boundary), the PDE (\ref{mfpt:ssp0}) is
radially symmetric and, from its readily obtained explicit solution,
an exact formula for the flux $J$ is \cite{Collins49}
\begin{equation}  \label{eq:J_kappa}
J = \int\limits_{\pa} (-D\partial_n {\mathcal U}) \, ds
= \frac{J_{\rm Smol}}{1 + D/(\K R)} \,,
\end{equation}
where $J_{\rm Smol} = 4\pi D R\, {\mathcal U}_{\rm inf}$ is the
Smoluchowski's diffusive flux onto a perfectly absorbing sphere
\cite{Smoluchowski1918}.   In order to quantify the effect of
a heterogeneous distribution of reactive patches, it is convenient to
determine the effective reactivity $\K_{\rm eff}$ of an equivalent
homogeneous sphere that produces the flux $J$ from (\ref{eq:J_kappa}).
In this way, by equating the fluxes in (\ref{eq:j_origin}) and
(\ref{eq:J_kappa}), we identify $\K_{\rm eff}$ in terms of
${\mathcal C}_{\rm T}$ as
\begin{equation}  \label{eq:kappa_J}
 \keff = \frac{R}{D} \K_{\rm eff} = \frac{1}{J_{\rm Smol}/J - 1} =
  \frac{1}{R/{\mathcal C}_{\rm T} - 1}\,.
\end{equation}

Berg and Purcell \cite{BergPurcell1977} pioneered the study of the
flux onto a spherical target that contains $N$ identical small
circular perfectly reactive boundary patches of radius $\eps R$. Based
on physical insight, they derived the following approximation for the
flux:
\begin{equation}  \label{eq:BP}
J_{\rm BP} = \frac{J_{\rm Smol}}{1 + \pi/(\eps N)}  \,.
\end{equation}
With this approximation, the effective reactivity of small patches
becomes ${R\K_{\rm BP}/D} = \eps N/\pi$, which can also be written in
terms of the patch area coverage fraction $f$ on the surface of the
sphere as (see also \cite{Shoup82})
\begin{equation}\label{eq:kapBP}
  \K_{\rm BP} = \frac{4 f D}{\pi \eps R} \,, \qquad
  \mbox{where} \qquad f\equiv\frac{N \pi (\eps R)^2}{4\pi R^2} =
  \frac{\eps^2 N}{4}
  \,.
\end{equation}
Qualitatively, this result is equivalent to the reactivity of $N$
one-sided disks of radius $\eps R$ in $\R^3$, as if they were trapping
diffusing particles independently of each other (in fact, as
$2 \eps R/\pi$ is the capacitance of a two-sided disk, half of it
corresponds to a one-sided disk).  The seminal approximation
(\ref{eq:BP}) ignores the diffusion screening
\cite{Felici03,Grebenkov05} or diffusive interactions \cite{Traytak92}
between patches that compete for the diffusing particles. This
competition is expected to reduce the fluxes to individual disks and
can thus significantly diminish the effective reactivity.  In
addition, the approximation (\ref{eq:kapBP}) ignores the curvature of
the spherical boundary that may also be relevant when the patches are
not too small.  Despite these limitations, the seminal work by Berg
and Purcell \cite{BergPurcell1977} stimulated the development of
asymptotic and homogenization methods in both the mathematical and
physical literature (see an overview in \cite{GrebSkvortsov}).

From a mathematical perspective, the method of matched asymptotic
expansions was used in \cite{Lindsay17} to provide a {\em first
principles} derivation of the effective reactivity for $N$ small
circular and perfectly reactive patches on the boundary of a sphere,
in the small patch area fraction limit $f\ll 1$.  This analysis
accounted for diffusive interactions between patches, as well the
curvature of the sphere. For uniformly distributed identical patches,
the leading-order term in the scaling law derived in \cite{Lindsay17}
was found to agree with the classical Berg-Purcell result
(\ref{eq:BP}). This asymptotic scaling law for the effective
reactivity was validated by comparing results with those computed from
a numerical PDE-based approach for the Dirichlet-Neumann BVP,
as developed in \cite{Lindsay17,Lindsay18a}, which effectively
resolved the edge singularity near each patch (see also
\cite{Eun2020}).  Moreover, a fast solver relying on an integral
equation re-formulation of the mixed Dirichlet-Neumann BVP for
a sphere with locally circular patches was developed in \cite{Kaye20}
to solve the PDE with up to $100,000$ perfectly reacting
patches. Other analytical and numerical studies with perfectly
reactive patches include \cite{Lindsay18a} and \cite{Lindsay18b} for
the boundary homogenization of periodic patterns of patches on a
semi-infinite plane (see also the references therein). More recently,
a hybrid analytical and numerical approach, based on a Kinetic Monte
Carlo algorithm, has been developed to study time-dependent diffusive
capture to an infinite plane and to a sphere \cite{Lindsay25}.

In contrast, there have been much fewer analytical or numerical
studies for the more realistic situation in (\ref{mfpt:ssp0}) of
partially reactive patches. In \cite{Zwanzig1991,Berez04}, a heuristic
interpolation formula was postulated, but without any mathematical
justification, for the effective reactivity of a generic smooth
surface that has identical small circular patches with a common
reactivity $\K$:
\begin{equation}\label{eq:kapheur}
  \K_{\rm heur} = \frac{f \K}{f\K +\K_{\rm BP}} \K_{\rm BP}\,,
\end{equation}
where $\K_{\rm BP}$ is given in (\ref{eq:kapBP}).  For identical
patches of radii $\eps R$ on the surface of sphere of radius $R$, the
heuristic approximation (\ref{eq:kapheur}) can be written
equivalently, by using (\ref{eq:kapBP}) for $\K_{\rm BP}$, as
\begin{equation}\label{eq:kapheur_2}  
  \frac{R}{D} \K_{\rm heur} = \frac{2 f}{\eps}
  {\mathcal A}\left( \frac{\K \eps R}{D} \right) ,
  \qquad   \mbox{where} \qquad  {\mathcal A}(\mu)\equiv
  \frac{{2\mu/\pi}}{\mu + {4/\pi}} \,.
\end{equation}
In \cite{Plunkett24} a {\em first principles} leading-order asymptotic
theory was developed to derive the effective reactivity for small
circular patches of radius $a$ on an infinite plane in $\R^3$ in the
limit of low, moderate, and large patch reactivity. With $a=\eps R$,
their leading-order analysis, resulting in Eq. (7.3) of
\cite{Plunkett24}, predicted that the effective reactivity can be
well-approximated, uniformly in $\K$, by (\ref{eq:kapheur_2}) with the
term ${4/\pi}$ in ${\mathcal A}(\mu)$ replaced by ${2/(\pi K_w)}$. The
value $K_w\approx 0.5854$ was estimated in \cite{Plunkett24} by using
Monte Carlo simulations to calculate the numerical solution to a
certain local PDE problem defined near a patch.  This numerical
value has now been corrected\footnote{As communicated to MJW by S.
Lawley, the Monte Carlo approximation should be $K_w\approx 0.5$ and
not $K_w\approx 0.5854$.} to $K_w\approx 0.5$, so that the results in
\cite{Plunkett24} do in fact provide the first theoretical
justification of the empirical approximation (\ref{eq:kapheur}) of
\cite{Zwanzig1991,Berez04}.

Related asymptotic studies for the mean first-passage time and
splitting probabilities for the narrow capture of a Brownian particle
in a bounded 3-D domain with small surface patches of finite
reactivity include \cite{Guerin23,Cengiz24}, and our companion paper
\cite{GrebWard25}. The analogous 2-D problem has been analyzed in
\cite{GrebWard25b}. Narrow capture in a 3-D bounded domain with
reflecting walls but with a collection of small spherical inclusions
with semipermeable interfaces has been analyzed in
\cite{Bressloff23f}.

For the sphere, our main goal is to extend the previous analysis in
\cite{Lindsay17}, which was restricted to perfectly reactive
patches, and the leading-order analysis in \cite{Plunkett24},  to
the more general case of partially reactive patches; in fact,
we aim at developing a three-term asymptotic formula for ${\mathcal
C}_{\rm T}$ in (\ref{mfpt:ssp0}) for patches of arbitrary shape
and any reactivity.  With this three-term asymptotic analysis we will
incorporate inter-patch interactions, as well as the effect of the
curvature of the sphere. By homogenizing our result for ${\mathcal
C}_{\rm T}$, we will derive a new scaling law for the effective
reactivity of the boundary of a sphere for identical partially
reactive patches. Results from the asymptotic theory will be confirmed
through simulations from a new Monte Carlo algorithm.  We remark that
the leading-order term in our asymptotic results is valid for an
arbitrary bounded domain ${\mathcal B}$ with a smooth boundary
$\pa$. For this extension, $\eps={L/R}$, where $2R$ is the diameter of
${\mathcal B}$, defined as the maximum Euclidean distance between any
two points on $\partial{\mathcal B}$.

As convenient for our analysis, we will first non-dimensionalize
(\ref{mfpt:ssp0}) in a sphere, by introducing the dimensionless
variables defined by
\begin{equation}\label{intro:scalings}
  \x = \frac{\X}{R} \,, \quad \Omega= \frac{{\mathcal B}}{R} \,, \quad
  \eps = \frac{L}{R}   \,, \quad a_i = \frac{L_i}{L} \,, \quad
  \kappa_i = \frac{L \K_i}{D} \,, \quad 
  u = \frac{{\mathcal U}}{{\mathcal U}_{\rm inf}}\,,
  \quad C_{\rm T}=\frac{{\mathcal C}_{\rm T}}{R} \,.
\end{equation}
In the region exterior to the unit sphere $\Omega$, we obtain that
(\ref{mfpt:ssp0}) transforms to
\bsub \label{mfpt:ssp}
\begin{align}
  \Delta_{\x} u & = 0 \,, \quad \x \in \R^3\backslash\Omega \,,
                  \label{mfpt:ssp_1}\\
  \eps\partial_{n} u + \kappa_i u & = 0\,, \quad \x \in \partial\Omega^{\eps}_i
                                    \,, \quad i=1,\ldots,N \,, \\
  \partial_{n} u & = 0 \,, \quad \x \in \partial \Omega_r
= \partial \Omega \backslash 
\overline{(\partial \Omega^{\eps}_1 \cup \ldots \cup \partial \Omega^{\eps}_N)}\,, \label{mfpt:ssp_2}\\
  u & \sim  1 - \frac{C_{\rm T}}{|\x|} +{\mathcal O}(|\x|^{-2})
  \,,  \quad \mbox{as} \quad |\x|\to \infty \,, \label{mfpt:ssp_3}
\end{align}
\esub where $\Delta_{\x}$ is the Laplacian in $\x$, and $\partial_n$
is the outward normal derivative to $\partial\Omega$, which points
into the unit sphere. In (\ref{mfpt:ssp}), each rescaled patch
$\partial\Omega^{\eps}_i = R^{-1} \partial{\mathcal B}_i$, of small
diameter ${\mathcal O}(\eps)$, is assumed to satisfy
$\partial\Omega^{\eps}_i\to\x_i\in \partial\Omega$ as $\eps\to 0$.
The patches are also assumed to be well-separated in the sense that
$\left|\x_i-\x_j\right|={\mathcal O}(1)$ for all $i\neq j$. 

In the limit $\eps\to 0$ of small patches, in \S
  \ref{mfpt_sec:expan} we will use strong localized perturbation
  theory, originating from \cite{Ward93}, to derive a three-term
  asymptotic expansion for the capacitance $C_{\rm T}$, valid for
  arbitrary $\kappa_i>0$, in which only the third-order term depends
  on the spatial arrangement $\lbrace{\x_1,\ldots,\x_N\rbrace}$ of the
  centers of the patches.  For the special case where
  $\kappa_i=\infty$ and when the patches are circular disks, such an
  asymptotic analysis has been performed in \cite{Lindsay17} by using
  the traditional spherical coordinates for a sphere. Our new
  framework in \S \ref{mfpt_sec:expan}, based on the geodesic normal
  coordinates introduced in \S \ref{mfpt_sec:prelim}, provides an
  alternative approach to that used in \cite{Lindsay17}.  These
  geodesic coordinates have been used previously in \cite{Tzou2020}
  and \cite{Gomez2020} to analyze localized pattern formation problems
  in reaction-diffusion systems. With this new, more efficient,
  coordinate system, we are readily able to extend the previous
  analysis of \cite{Lindsay17} to the case of finite $\kappa_i$ and to
  arbitrary patch shapes.  Our main result, summarized in Proposition
  \ref{mfpt_b:main_res} of \S \ref{mfpt_sec:expan}, involves the
  reactive capacitance $C_i(\kappa_i)$ and an additional monopole
  coefficient $E_{i}(\kappa_i)$, which are determined in terms of the
  solutions to certain local problems near each patch. We emphasize
  that up to three terms in the asymptoptic expansion our result in
  Proposition \ref{mfpt_b:main_res} is identical in form to that
  derived in Eq. (3.37) of \cite{Lindsay17} for circular patches of
  perfect reactivity provided that we simply replace the two monopole
  coefficients in Eq. (3.37) of \cite{Lindsay17} with $C_i(\kappa_i)$
  and $E_{i}(\kappa_i)$.

From the analysis in the companion paper \cite{GrebWard25}, we show
that $C_i(\kappa_i)$ and $E_i(\kappa_i)$ can be expressed for all
$\kappa_i>0$ via a Steklov eigenfunction expansion.  For circular
patches, an efficient numerical computation of the Steklov
eigenfunctions is readily available \cite{Grebenkov24}.  Moreover,
both $C_i(\kappa_i)$ and $E_i(\kappa_i)$ can be well-approximated over
the full range $0<\kappa_i<\infty$ by some heuristic formulas that
have been benchmarked to full numerical results. In this way, we
obtain an explicit analytical approximation for $C_{\rm T}$ for any
$\kappa_i>0$. Preliminary results for $C_i(\kappa_i)$,
$E_i(\kappa_i)$, and the surface Neumann Green's function, which are
central to the analysis in \S \ref{mfpt_sec:expan}, are summarized in
\S \ref{mfpt_sec:prelim}.  Although some of these results were
previously derived in \cite{GrebWard25}, they are summarized here for
completeness and readability.

In \S \ref{sec:homog} we consider the homogenization limit of $N\gg 1$
identical patches that are uniformly distributed over the surface of
the sphere, but in the limit of small patch area fraction. For this
homogenization problem, we derive an explicit analytical scaling law
for both the effective capacitance $\ceff$ and the effective
reactivity $\keff$ of the structured target.  This scaling law can be
applied over the full range $\kappa>0$ by using our empirical
approximations for $C(\kappa)$ and $E(\kappa)$.  Substitution of the
empirical approximation for $C(\kappa)$ into the leading-order term of
our asymptotic theory leads to the heuristic result
(\ref{eq:kapheur_2}) of \cite{Berez04}.  In \S \ref{sec:numer} we
present a new efficient Monte Carlo method to numerically compute the
capacitance $C_{\rm T}$ from the underlying PDE.  From numerical results
obtained by the Monte Carlo algorithm, we show in \S \ref{sec:comp}
that our explicit scaling laws for $\ceff$ and $\keff$, which are
based on our three-term asymptotic theory, are still accurate even
when there is only a moderately large number of patches, or when the
reactive patch area fraction is not exceedingly
small. Finally, in \S \ref{sec:discussion} we summarize
 our results and discuss a few additional problems
  that warrant further study.

\section{Preliminaries}\label{mfpt_sec:prelim}

In this section we derive some preliminary results that are central
for our asymptotic analysis in \S \ref{mfpt_sec:expan}.

\subsection{Geodesic normal coordinates}

The inner problem near each patch is more readily analyzed in terms of
geodesic normal coordinates, rather than the usual global spherical
coordinates used in \cite{Lindsay17}. More specifically, for 
each patch $\partial\Omega_i^\eps$ on the unit sphere $\Omega$, we
introduce the geodesic normal coordinates
$\bxi=(\xi_1,\xi_2,\xi_3)^T\in \left({-\pi/2},{\pi/2}\right) \times
\left(-\pi,\pi\right)\times [0,\infty]$ in $\R^3\backslash\Omega$ so
that $\bxi=0$ corresponds to $\x_i\in\partial\Omega$, while $\xi_3>0$
corresponds to the exterior of $\Omega$.  In these coordinates,
$\xi_1$ can be viewed as the polar angle of a spherical coordinate
system centered at $\x_i$ on the sphere, but defined on the range
$\xi_1\in\left({-\pi/2},{\pi/2}\right)$ that avoids the usual
coordinate singularity of spherical coordinates at the north pole.
The curves obtained by setting $\xi_3=0$ and fixing either $\xi_1=0$
or $\xi_2=0$ are geodesics on $\partial\Omega$ that pass through
$\x_i$.

\begin{figure}[htbp]
  \centering
  \includegraphics[width=0.35\textwidth]{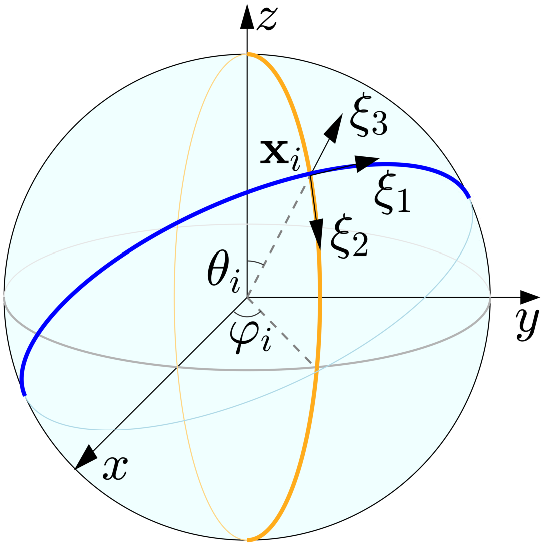} 
    \caption{Geodesic normal coordinates $(\xi_1,\xi_2,\xi_3)^T$ at
      $\x_i\in \partial\Omega$, with the geodesics (orange and blue
      curves) indicated. }
    \label{fig:geodesic}
\end{figure}

From the global transformation $\x=\x(\bxi)$ between cartesian and
geodesic coordinates, as given explicitly in (\ref{app_g:global}) of
Appendix \ref{app_g:geod}, we derive an exact expression for the
Laplacian of a generic function ${\mathcal V}(\bxi)\equiv
u\left(\x(\bxi)\right)$. Then, in terms of the local (or inner)
variables, $\y=(y_1,y_2,y_3)^T$, defined by
\begin{equation}\label{mfpt:innvar}
  \xi_1=\eps y_1 \,, \qquad \xi_2=\eps y_2 \,, \qquad \xi_3 =\eps y_3\,,
\end{equation}
in Appendix \ref{app_g:geod} we show, for $\eps\to 0$ and for
$V(\y)={\mathcal V} (\eps \y)$, that
\begin{equation}\label{mfpt:local}
  \Delta_{\x} u = \eps^{-2} \Delta_{\y} V + \eps^{-1} \left[-2 y_3 \left(
      V_{y_1 y_1} + V_{y_2 y_2}\right) + 2 V_{y_3} \right] +
  {\mathcal O}(1)\,,
\end{equation}
with $\Delta_{\y}V \equiv V_{y_1 y_2} + V_{y_2 y_2} + V_{y_3
  y_3}$. This result is essential to our local analysis in \S
\ref{mfpt_sec:expan}.

\subsection{Reactive capacitance}

Defined on the tangent plane to the sphere centered at $\x_i$, the
leading-order term in our inner expansion near $\x=\x_i$ will involve
the solution $w_i = w_{i}(\y;\kappa_i)$ to
\bsub \label{mfpt:wc}
\begin{align}
    \Delta_{\y} w_{i} &=0 \,, \quad \y \in \R_{+}^{3} \,, \label{mfpt:wc_1}\\
    -\partial_{y_3} w_{i} + \kappa_i w_{i} &=\kappa_i \,, \quad y_3=0 \,,\,
    (y_1,y_2)\in \PT_i\,,  \label{mfpt:wc_2}\\
    \partial_{y_3} w_{i} &=0 \,, \quad y_3=0 \,,\, (y_1,y_2)\notin \PT_i
    \,, \label{mfpt:wc_3}\\
  w_{i}&\sim \frac{C_{i}(\kappa_i)}{|\y|} +
                {  \frac{\DT_i(\kappa_i) {\bf \cdot} \y}{|\y|^3}}
         + \cdots\,,  \quad \mbox{as}\quad
    |\y|\to \infty \,. \label{mfpt:wc_4}
\end{align}
\esub Here
$\R_{+}^{3}\equiv\lbrace{\y=(y_1,y_2,y_3) \, \vert \, \, y_3\geq 0 \,,
  \, -\infty<y_1,y_2<\infty \rbrace}$ is the upper half-space, and
$\PT_i \asymp \eps^{-1}\partial\Omega^{\eps}_i$ is the flattened Robin
patch on the horizontal plane $y_3=0$.  In (\ref{mfpt:wc_4}), the
dipole vector $\DT_i=\DT_i(\kappa_i)$ must have the form
$\DT_i=(p_{1i},p_{2i},0)^T$ to ensure that the far-field behavior
(\ref{mfpt:wc_4}) satisfies (\ref{mfpt:wc_3}). When $\Gamma_i$ is
symmetric in $y_1$ and $y_2$, such as is the case when $\Gamma_i$ is a
disk, it follows from symmetry that $p_{i1}=p_{i2}=0$, so that the
dipole term in the far-field (\ref{mfpt:wc_4}) is absent.

We refer to the monopole coefficient $C_i(\kappa_i)$ in
(\ref{mfpt:wc_4}) as the {\em reactive capacitance}. Although by
applying the divergence theorem we can readily  express
$C_i(\kappa_i)$ in terms of the {\em charge density} $q_i$, 
\begin{equation}\label{mfpt:wc_charge}
  C_i(\kappa_i) = \frac{1}{\pi} \int_{\PT_i} q_{i}(y_1,y_2; \kappa_i) \,
  dy_1 dy_2 , \qquad \mbox{where} \qquad
  q_i(y_1,y_2;\kappa_i)\equiv -\frac{1}{2} \partial_{y_3} w_{i}\vert_{y_3=0}\,,
\end{equation}
we emphasize that there is no {\it explicit} solution to
(\ref{mfpt:wc}) for arbitrary $\kappa_i>0$.  As a result,
$C_i(\kappa_i)$ must in general be computed numerically. In Appendix
\ref{app:Cmu} we summarize some results of Appendix D of
\cite{GrebWard25} that determined an easily computable spectral
representation for $C_i(\kappa_i)$ in terms of a suitable Steklov
eigenvalue problem. More specifically, this spectral representation is
\begin{equation}  \label{eq:Cmu_def0}
  C_i(\kappa_i) = \frac{\kappa_i}{2\pi} \sum\limits_{k=0}^\infty
\frac{\mu_{ki} d_{ki}^2}{\mu_{ki} + \kappa_i} \,,
\end{equation}
where the eigenvalues $\mu_{ki} > 0$ and the spectral weights $d_{ki}
\ne 0$ of the Steklov problem (\ref{eq:Psi_def}) can be computed
numerically for a specified patch shape $\Gamma_i$. Although the
spectrum of (\ref{eq:Psi_def}) must be computed numerically, the
functional form of $C_{i}(\kappa_i)$ and its dependence on reactivity
is universal.  In the special case where the patches are all of the
same shape, but with a variable size, we have the scaling law
\begin{equation} \label{eq:Cmu_def0_1}
  C_i(\kappa_i)=a_i {\mathcal C}(\kappa_i a_i) \,, \quad i=1\,,\ldots\,,N,
\end{equation}
where ${\mathcal C}(\mu)$ is the reactive capacitance of the rescaled
common patch shape $\Gamma_c\equiv {\Gamma_i/a_i}$, which needs to be
computed only once. In particular, this scaling result is applicable
to the case where the patches are all circular disks of different
radii. For a circular patch of unit radius, the first eight Steklov
eigenvalues and weights, corresponding to eigenfunctions that are
axially symmetric on the patch, are given in Table
\ref{table:muk_disk}.

\subsubsection*{Circular patch}

We recall that when $\Gamma_i$ is a disk of radius $a_i$, the limiting
problem with $\kappa_i=\infty$ is the classical problem for the
capacitance of a flat disk (cf.~\cite{Jackson}), whose solution
$w_i(\y;\infty)$ is (see page 38 of \cite{Fabrikant89})
\bsub \label{mfpt:wcinf}
\begin{equation}
  w_i(\y;\infty) \equiv \frac{2}{\pi}
  \sin^{-1}\left(\frac{a_i}{B(y_3,\rho_0)} \right) \,,
    \label{mfpt:wcinf_1}
\end{equation}
where $\rho_0 \equiv (y_1^2 + y_2^2)^{1/2}$ and 
\begin{equation}
  B(y_3,\rho_0) \equiv  \frac{1}{2} \left(  \left[ (\rho_0 + a_i)^2 + y_3^2
  \right]^{1/2}  +   \left[ (\rho_0 - a_i)^2 + y_3^2 \right]^{1/2} \right) \,.
  \label{mfpt:wcinf_2}
\end{equation}
The far-field behavior is
$w_i(\y;\infty)\sim {C_{i}(\infty)/|\y|}+{\mathcal O}(|\y|^{-3})$, where
$C_i(\infty)={2a_i/\pi}$. Moreover, from (\ref{mfpt:wcinf_1}), and by
using the radial symmetry, the charge density is
\begin{equation}\label{mfpt:wc_q}
 q_i(y_1,y_2;\infty)=q_i(\rho_0;\infty)\equiv
  -\frac{1}{2} \partial_{y_3} w_i(\y;\infty)\vert_{y_3=0}=
  \frac{1}{\pi \sqrt{a_i^2-\rho_0^2}} \,, \quad 0\leq \rho_0\leq a_i \,.
\end{equation}
\esub
As a result of the edge singularity in (\ref{mfpt:wc_q}) at
$\rho=a_i$, it can be shown from the analysis in \cite{Guerin23} (see
also Appendix D.3 of \cite{GrebWard25}) that the difference
$C_i(\kappa_i)-C_{i}(\infty)$ is not analytic for $\kappa_i\gg 1$. More
specifically, $C_i(\kappa_i)$ has the refined far-field behavior
\begin{equation} \label{eq:Cmu_asympt}
  C_i(\kappa_i) \sim \frac{2a_i}{\pi} - 2 \frac{\left[\log(a_i\kappa_i) +
      \log{2}+\gamma_e+1\right]}{\pi^2\kappa_i} \,, \quad \mbox{as}\quad
  \kappa_i\to\infty \,,
\end{equation}
where $\gamma_e \approx 0.5772\ldots$ is Euler's constant.

\begin{figure}[htbp]
  \centering
     \begin{subfigure}[b]{0.49\textwidth}  
      \includegraphics[width =\textwidth]{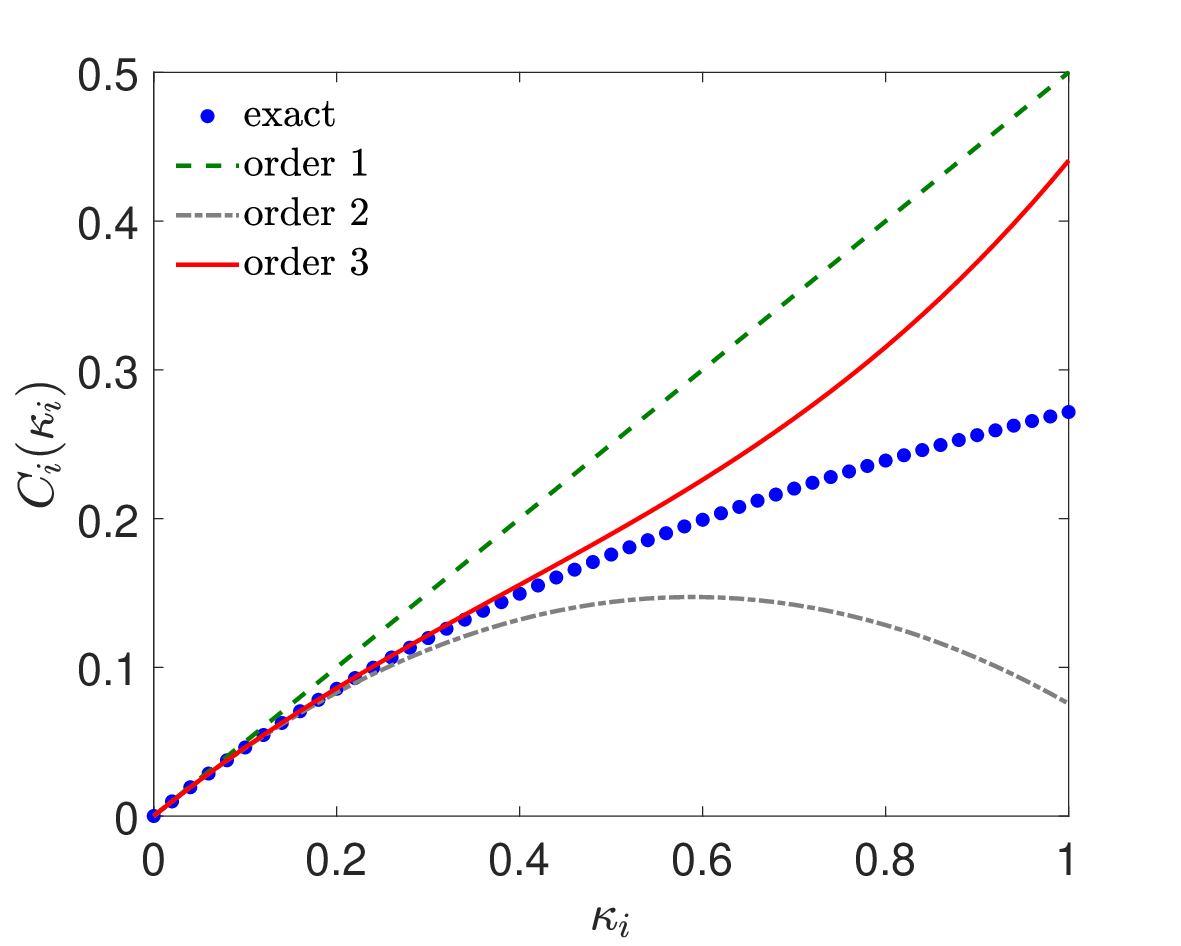} 
        \caption{Small $\kappa_i$ comparison}
        \label{fig:Cmu_small}
    \end{subfigure}
    \begin{subfigure}[b]{0.49\textwidth}
      \includegraphics[width=\textwidth]{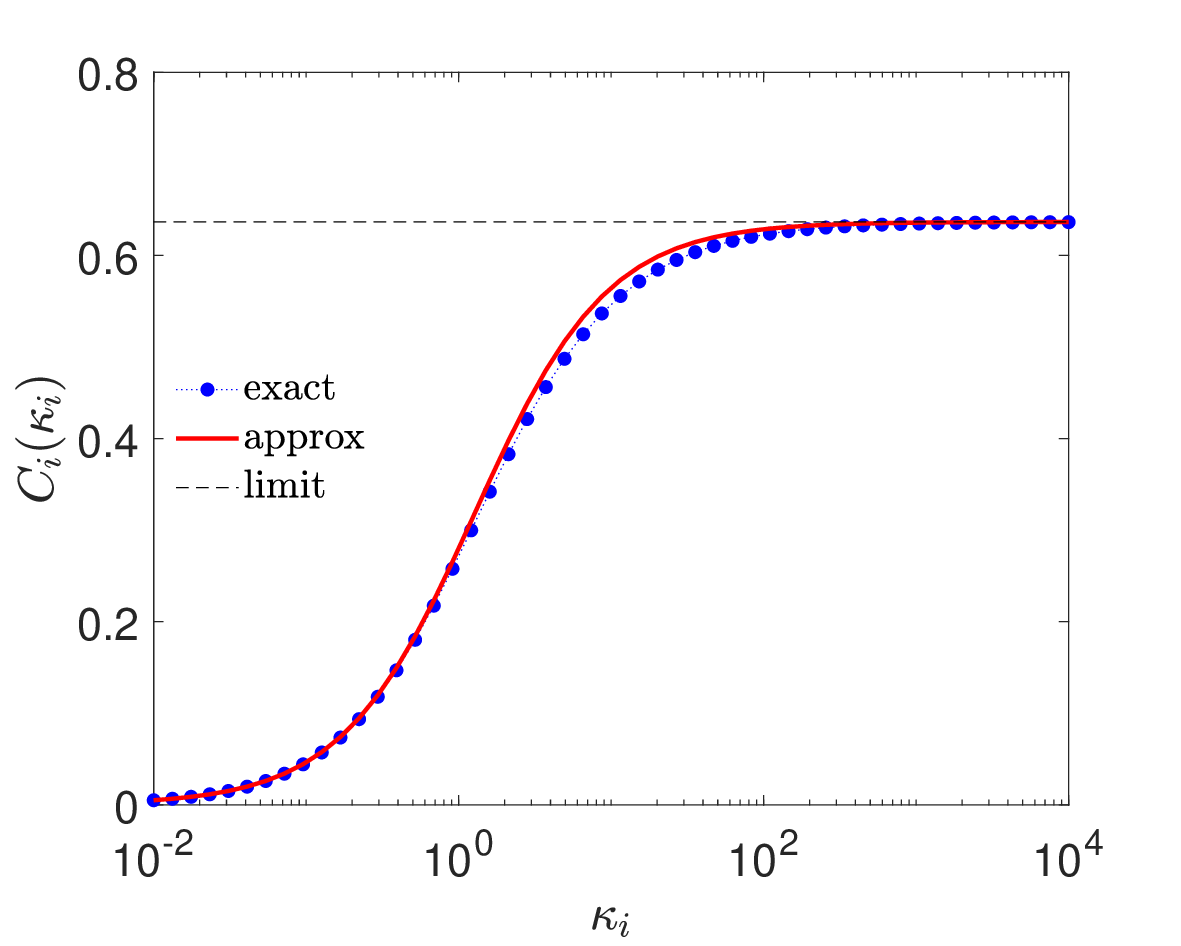} 
        \caption{Sigmoidal approximation} 
        \label{fig:Cmu_approx_new}
    \end{subfigure}
\caption{
The reactive capacitance $C_i(\kappa_i)$ for the circular patch
$\PT_i$ of unit radius ($a_i = 1$).  {\bf (a):} A comparison of
$C_i(\kappa_i)$ numerically computed from (\ref{eq:Cmu_def0}), with
the one-, two-, and three-term approximations obtained from
(\ref{eq:Cmu_Taylor}) and (\ref{eq:cn_exact}), valid for $\kappa_i \ll
1$.  {\bf (b):} The sigmoidal approximation (\ref{mfpt:sigmoidal_2})
provides a decent approximation of the numerical result for
$C_i(\kappa_i)$ on the full range $\kappa_i >0$.}
\label{fig:capprox}
\end{figure}

In Fig.~\ref{fig:Cmu_small} we illustrate that a three-term
Taylor expansion (\ref{eq:Cmu_Taylor}) of $C_i(\kappa_i)$ at
small $\kappa_i$, with the coefficients in (\ref{eq:cn_exact}),
provides a rather close approximation to $C_i(\kappa_i)$ on the range
$0<\kappa_i<0.45$. Finally, in Fig.~\ref{fig:Cmu_approx_new} we show
that the heuristic sigmoidal approximation, given by
\begin{equation}\label{mfpt:sigmoidal_2}
  C_{i}(\kappa_i)\approx C_i^{\rm app}(\kappa_i) =
  a_i {\mathcal C}^{\rm app}(a_i\kappa_i)\,,  \qquad
  \mbox{where} \qquad {\mathcal C}^{\rm app}(\mu)= \frac{{2\mu/\pi}}
  {\mu + {4/\pi}} \,, 
\end{equation}
is within $4\%$, over the entire range $\kappa_i>0$, of the spectral
expansion result (\ref{eq:Cmu_def0}) applied to a circular patch.
These results are summarized in the following lemma of
\cite{GrebWard25}.

\begin{lemma}(Lemma 2.1 of \cite{GrebWard25}) \label{lemma:Cj_kappa} 
When $\Gamma_i$ is the disk $y_1^2+y_2^{2}\leq a_i^2$, $C_i(\kappa_i)$
has the limiting asymptotics
\bsub  \label{mfpt:cj_small}
  \begin{align}
    C_i(\kappa_i) &\sim C_i(\infty) + {\mathcal O}\left(
  \frac{\log\kappa_i}{\kappa_i}\right) \,,\quad \mbox{as}\quad
                    \kappa_i\to\infty\,, \quad \mbox{with} \quad
  C_i(\infty)=\frac{2a_i}{\pi} \,, \label{mfpt:cj_large_b}  \\  
  C_i(\kappa_i) &  \sim a_i \biggl[c_{1i} \kappa_i a_i - c_{2i} (\kappa_i a_i)^2 
     + c_{3i} (\kappa_i a_i)^3 + {\mathcal O}((\kappa_ia_i)^4)\biggr] \,,
                  \quad \mbox{as}
  \quad \kappa_i\to 0 \,,  \label{mfpt:cj_small_b}
  \end{align}
  \esub where $c_{1i}=0.5$, $c_{2i} \approx 0.4241$ and
  $c_{3i} \approx 0.3651$ (see (\ref{eq:cn_exact}) of Appendix
  \ref{app:Cmu}), are independent of the patch radius $a_i$. The
  sigmoidal approximation (\ref{mfpt:sigmoidal_2}) is exactly
  consistent with only the leading-order coefficient $c_{1i}$. The
  other Taylor coefficients $c_{2i}$ and $c_{3i}$ are numerically
  comparable, but distinct, from the higher Taylor coefficients of
  (\ref{mfpt:sigmoidal_2}).
\end{lemma}

\subsection{Monopole From a Higher-Order Inner Solution}
\label{prel:high}

In our asymptotic analysis of (\ref{mfpt:ssp}) in \S
\ref{mfpt_sec:expan} we show in Appendix \ref{app_h:inn2} that we have to
calculate the monopole coefficient $E_i=E_{i}(\kappa_i)$ defined by
the solution $\Phi_{2hi}$ to the inner problem
\bsub\label{mfpt:inn2_probh}
\begin{align}
\Delta_{\y} \Phi_{2hi} &= 0 \,, \quad \y \in \R_{+}^{3} \,,\label{mfpt:inn2_h1}\\
  -\partial_{y_3} \Phi_{2hi} + \kappa_i \Phi_{2hi} &= \kappa_i {\mathcal F}_{i}
                                                 \,, \quad y_3=0 \,,\,
    (y_1,y_2)\in \PT_i\,,  \label{mfpt:inn2_h2}\\
    \partial_{y_3} \Phi_{2hi} &=0 \,, \quad y_3=0 \,,\, (y_1,y_2)\notin \PT_i
    \,, \label{mfpt:inn2_h3}\\
    \Phi_{2hi} & \sim -\frac{E_i}{|\y|}\,,  \quad
    \mbox{as} \quad |\y| \to \infty \,. \label{mfpt:inn2_h4}
\end{align}
\esub In (\ref{mfpt:inn2_h2}), the inhomogeneous term
${\mathcal F}_i={\mathcal F}_{i}(y_1,y_2;\kappa_i)$ is the unique
solution, defined in terms of $C_i=C_i(\kappa_i)$ and
  $q_i(y_1,y_2;\kappa_i)$ as related by (\ref{mfpt:wc_charge}), to
the 2-D problem
\bsub\label{mfpt:fprob}
\begin{gather}  
  {\mathcal F}_{i,y_1 y_1} + {\mathcal F}_{i,y_2 y_2}=q_i(y_1,y_2;\kappa_i)
  I_{\PT_i}(y_1,y_2) \,,
    \quad I_{\PT_i} \equiv \left\{\begin{array}{ll}
        1 \,, & (y_1,y_2) \in \PT_i \\
    0 \,, & (y_1,y_2) \notin \PT_i  \end{array}\right. \,,
  \label{mfpt:fprob_1}\\
  {\mathcal F}_{i} \sim \frac{C_i}{2}\log\rho_0 + o(1)\,,  \quad
  \mbox{as} \quad \rho_0\equiv (y_1^2+y_2^2)^{1/2}\to \infty
                            \,. \label{mfpt:fprob_2}
\end{gather}
\esub It is the $o(1)$ condition in the far-field
  (\ref{mfpt:fprob_2}) which ensures that ${\mathcal F}_i$ is unique.

For an arbitrary patch shape, in Appendix \ref{app_h:inn2} we show in
(\ref{app_h:eval}) that $E_i(\kappa_i)$ can be determined up to a
quadrature. The following result, established in \cite{GrebWard25},
and discussed in Appendix \ref{app_h:inn2}, more fully characterizes
$E_i$ when $\PT_i$ is a disk:

\begin{lemma}(Lemma 2.2 of
  \cite{GrebWard25}) \label{lemma:Ej_kappa}
  When the Robin patch $\Gamma_i$ is the disk $y_1^2+y_2^{2}\leq a_i^2$, we
  have
  \begin{equation}\label{mfpt:Ej_all}  
  E_i =  E_i(\kappa_i) = -\frac{\log{a_i}}{2} [C_i(\kappa_i)]^2 + a_i^2 \,
  {\mathcal E}_i(\kappa_i a_i) \,,
\end{equation}
where $C_i=C_i(\kappa)$ is the reactive capacitance for a disk given by
\begin{equation}\label{mfpt:Ej_part}
  C_i= 2 \int_{0}^{a_i} q_i(\rho_0;\kappa_i) \rho_0 \, d\rho_0 \,, \qquad
  q_i(\rho_0;\kappa_i)=-\frac{1}{2} w_{i, y_3}\vert_{y_3=0} \,,
\end{equation}
and ${\mathcal E}_i(\mu)$ is defined by
\begin{equation} 
{\mathcal E}_i(\mu) \equiv 2 \int_{0}^1 \frac{1}{\rho_0}
\left(\int_{0}^{\rho_0}  a_i q_i(\eta a_i;\mu/a_i) \, \eta \,
  d\eta\right)^2 \, d\rho_0 \,.
\end{equation}
The limiting asymptotic behavior of $E_i(\kappa_i)$ is
\bsub \label{mfpt:Ej_asy}
  \begin{align}
  E_{i} &\sim E_{i}(\infty)\equiv -\frac{2 a_i^2}{\pi^2}\left(
    \log{a_i} + \log{4} - \frac{3}{2} \right)\,,   \quad \mbox{as}\quad
  \kappa_i\to \infty\,, \label{mfpt:EJ_asy_1}\\
  E_i & \sim \frac{\kappa_i^2 a_i^4}{8} \left( \frac{1}{4}-\log{a_i}\right)\,,
  \quad \mbox{as} \quad \kappa_i\to 0 \,. \label{esmall}
\end{align}
\esub
\end{lemma}

\begin{figure}
\begin{center}
\includegraphics[width=88mm]{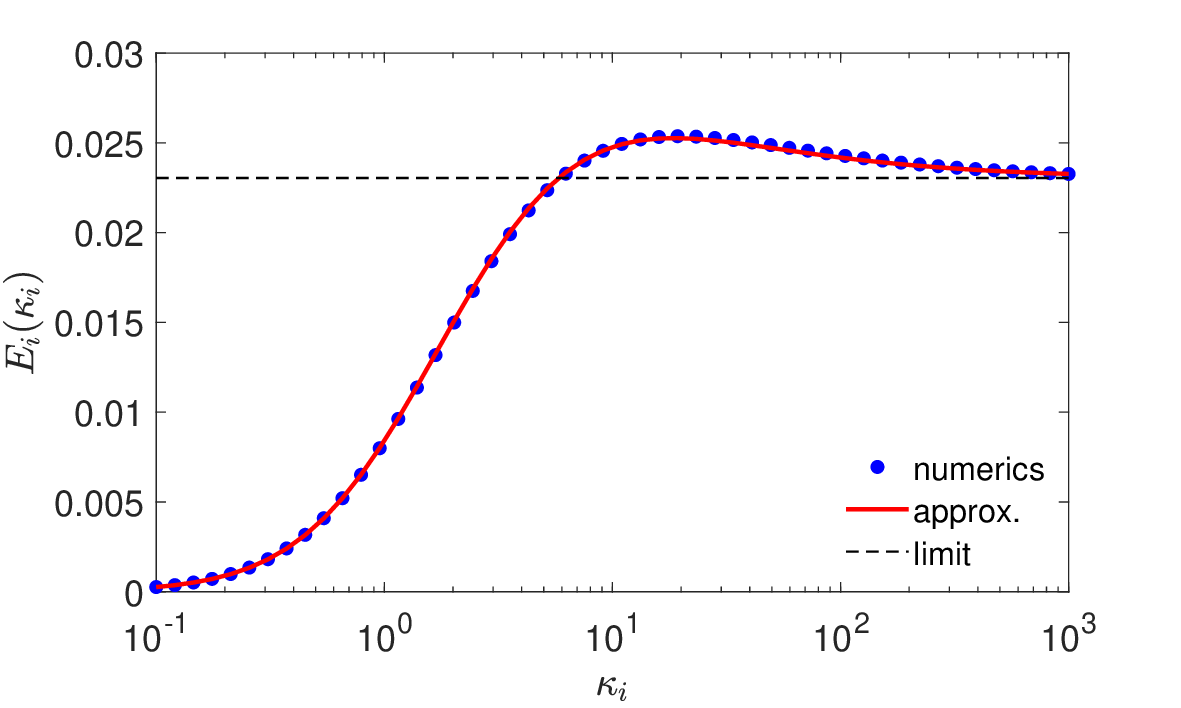} 
\end{center}
\caption{For a circular patch of radius $a_i=1$, the heuristic
  approximation ${\mathcal E}^{\rm app}(\kappa_i)$ (solid curve) from
  (\ref{mfpt:Eapp}) is compared with
  $E_i(\kappa_i) = {\mathcal E}_i(\kappa_i)$ (filled circles) given by
  (\ref{mfpt:Ej_all}) and computed via the numerical approach
  described in Appendix D of \cite{GrebWard25}. The dashed horizontal
  line is the asymptotic limiting value $(3 - 4\log 2)/\pi^2$
  consistent with (\ref{mfpt:EJ_asy_1}).}
\label{fig:Ekappa}
\end{figure}

Since there is no analytical formula for $E_i(\kappa_i)$ for arbitrary
$\kappa_i>0$ when $\Gamma_i$ is a disk, in Appendix D of
\cite{GrebWard25} we showed how it can be computed numerically to high
precision by expanding the charge density $q_i$ in terms of Steklov
eigenfunctions (see \cite{GrebWard25} for details). Moreover, labeling
$\mu\equiv a_i\kappa_i$ and with ${\mathcal C}^{\rm app}(\mu)$ as
given by the sigmoidal approximation (\ref{mfpt:sigmoidal_2}), in
Appendix D of \cite{GrebWard25} we showed that the heuristic
approximation \bsub \label{mfpt:E_heur}
\begin{equation}\label{mfpt:E_heur_1}
  E_{i}(\kappa_i)\approx E_{i}^{\rm app} = -\frac{a_i^2\log{a_i}}{2}
  \left[{\mathcal C}^{\rm app}(\kappa_i a_i)\right]^2 + a_i^2
  {\mathcal E}^{\rm app}(\kappa_i a_i)\,,
\end{equation}
where 
\begin{equation}  \label{mfpt:Eapp}
{\mathcal E}^{\rm app}(\mu) \equiv \left[{\mathcal C}^{\rm app}(\mu)\right]^2 
\biggl(\frac34 - \log{2} + \frac{1}{\frac{1}{\log{2} - 5/8} + 5.17 \,
  \mu^{0.81}} \biggr)\,,
\end{equation}
\esub agrees with the corresponding numerical result, with a maximal
relative error of $0.7\%$ over the entire range of $\mu >0$ (see
Fig.~\ref{fig:Ekappa}). Moreover, the established limits from
(\ref{mfpt:Ej_asy}) as $\mu\to 0$ and $\mu \to \infty$ are satisfied by
(\ref{mfpt:E_heur}).

A theoretical analytical explanation for the non-monotonic
  behavior of $E_i(\kappa_i)$ versus $\kappa_i$ observed in
  Fig.~\ref{fig:Ekappa} is not clear. It may be attributable to the
  non-monotonic ${\mathcal O}\left(\kappa_i^{-1}\log{\kappa_i}\right)$
  correction term in (\ref{mfpt:cj_large_b}) for $C_i$ when
  $\kappa_i\gg 1$ or the intricate behavior of the charge density
  $q_i(r;\kappa_i)$, defined in (\ref{mfpt:Ej_part}), as shown in
  Fig. E.1 of \cite{GrebWard25}. As derived in detail in
  \cite{Guerin23}, for $\kappa_i$ large but finite, this charge
  density is finite at $r=1$, but develops a boundary layer structure
  near $r=1$ as $\kappa_i$ is increased.

\subsection{The Exterior Surface Neumann Green's function}\label{prel:gs}

In our asymptotic analysis in \S \ref{mfpt_sec:expan},
the asymptotic expansion of the outer solution is represented in terms of
the {\em exterior} surface Neumann Green's function $G_s(\x;\x_i)$ for
the unit sphere $\Omega$, defined as the unique solution to
\begin{equation}
\Delta G_s = 0 \,, \quad  r>1\,; \quad
-\partial_r G_s =  \delta(\x-\x_i)\,, \quad r=1\,; \quad
G_s \sim \frac{1}{4\pi|\x|} \,, \quad \mbox{as} \quad |\x|\to \infty\,,
\label{mfpt:sph}
\end{equation}
where $r=|\x|$ and $|\x_i|=1$.  The exact solution to (\ref{mfpt:sph})
is (see \cite{Lindsay17}) 
\begin{equation}\label{mfpt:gs_exact}
G_s(\x;\x_i) = \frac{1}{2 \pi \left|\x-\x_i\right|} -
\frac{1}{4\pi } \log\left( 1 + \frac{2}{|\x-\x_i|+|\x|-1} \right)\,.
\end{equation}

To determine the local behavior of $G_{s}(\x;\x_i)$ as $\x\to\x_i$ in
terms of the local geodesic coordinates $\y$, we use (\ref{app_g:loc})
to estimate $|\x-\x_i|$. In this way, we obtain that
\begin{equation*}
  G_s \sim \frac{1}{2\pi\eps|\y|} \left( 1 - \frac{\eps y_3}{2|\y|^2}
    (y_1^2+y_2^2)\right) -\frac{1}{4\pi} \log\left(1 + \frac{2}{\eps (
      |\y| + y_3)}\right)\,,
\end{equation*}
which can be simplified asymptotically to
\begin{equation}\label{mfpt:gs_locm}
    G_s\sim \frac{1}{2\pi \eps |\y|} + \frac{1}{4\pi}\log\left(\frac{\eps}{2}
    \right) - \frac{y_3(y_1^2+y_2^2)}{4\pi |\y|^3} + \frac{1}{4\pi}
    \log(|\y| + y_3) + o(1)\,.
\end{equation}

\section{Analysis for the Capacitance $C_{\rm T}$}\label{mfpt_sec:expan}

We now use the method of matched asymptotic expansions to construct
solutions to (\ref{mfpt:ssp}) in the limit $\eps \to 0$. For our
analysis it is convenient to introduce $U$ by $u=-C_{\rm T} U$, so that
from (\ref{mfpt:ssp}) $U$ satisfies
\bsub \label{bp:ssp}
\begin{align}
  \Delta_{\x} U & = 0 \,, \quad \x \in \R^3\backslash\Omega \,,
                  \label{bp:ssp_1}\\
  \eps\partial_{n} U + \kappa_i U & = 0\,, \quad \x \in \partial\Omega^{\eps}_i
                                    \,, \quad i=1,\ldots,N \,, \\
  \partial_{n} U & = 0 \,, \quad \x \in \partial \Omega_r  \,, \label{bp:ssp_2}\\
  U & \sim  -\frac{1}{C_{\rm T}} + \frac{1}{|\x|} + \frac{\DT {\bf \cdot} \x}
 {C_{\rm T}|\x|^3} + \cdots\,, \quad \mbox{as}\quad |\x|\to \infty \,.\label{bp:ssp_3}
\end{align}
In this way, we can write the far-field condition (\ref{bp:ssp_3})
as the following flux condition over the boundary $\partial\Omega_R$ of a
large sphere of radius $R$ centered at $\x=0$:
\begin{equation}\label{bp:flux}
  \lim_{R\to\infty} \int_{\partial\Omega_R} \partial_{r} U \vert_{r=R} \, ds = - 4\pi\,.
\end{equation}
\esub

In the outer region away from the Robin patches we expand the outer
solution as
\begin{equation}
    U \sim \eps^{-1} U_0 + U_1 + \eps \log\left( \frac{\eps}{2} \right)
   U_2 + \eps U_3 + \cdots \,, \label{mfpt_b:outex}
\end{equation}
where $U_0$ is a constant to be determined, and where $U_k$ for $k\geq 1$
satisfies
\begin{equation}  \label{mfpt_b:Uk}
\begin{split}
&  \Delta_{\x} U_k = 0 \,, \quad \x \in \R^{3}\backslash \Omega \,; 
  \quad \partial_n U_k = 0 \,, \quad \x\in \partial\Omega\backslash
  \lbrace{\x_1,\ldots,\x_N\rbrace} \,, \\
  &  \lim_{R\to\infty} \int_{\partial\Omega_R} \partial_n U_k \vert_{r=R}\,
  ds = - 4\pi \delta_{k1}\,.
\end{split}
\end{equation}
Here $\delta_{k1}$ is the Kronecker symbol.  Our analysis will provide
singularity behaviors for each $U_{k}$ as $\x\to \x_i$, for
$i=1,\ldots,N$.

In the inner region near the $i$-th Robin patch we introduce the
local geodesic coordinates (\ref{mfpt:innvar}) and we expand each
inner solution as
\begin{equation}
 U \sim \eps^{-1} V_{0i} + \log\left( \frac{\eps}{2} \right) V_{1i} +
  V_{2i}  + \ldots \,. \label{mfpt_b:innex}
\end{equation}
Upon substituting (\ref{mfpt_b:innex}) into (\ref{mfpt:local}), we obtain
that $V_{ki}$ for $k=0,1,2$ satisfies
\bsub \label{mfpt_b:Vk}
\begin{align}
  \Delta_{\y} V_{ki} &= \delta_{k2} \left(- 2y_{3} V_{0i,y_3 y_3} - 2 V_{0i,y_3}
                       \right) \,, \quad
   \y \in \R_{+}^{3} \,, \label{mfpt_b:Vk_1}\\
   -\partial_{y_3} V_{ki} + \kappa_i V_{ki} &=0 \,, \quad y_3=0 \,,\,
    (y_1,y_2)\in \PT_i\,,  \label{mfpt_b:Vk_2}\\
    \partial_{y_3} V_{ki} &=0 \,, \quad y_3=0 \,,\, (y_1,y_2)\notin \PT_i
    \,. \label{mfpt_b:Vk_3}
\end{align}
\esub  

The leading-order matching condition is that $V_{0i}\sim U_0$ as
$|\y|\to\infty$ for each $i=1,\ldots,N$. As a result, we write the
leading-order inner solution as
\begin{equation}\label{mfpt_b:v0sol}
    V_{0i} = U_0 \left( 1 - w_{i} \right) \,,
\end{equation}
where $w_{i}=w_{i}(\y;\kappa_i)$ is the solution to (\ref{mfpt:wc}).
The asymptotic matching condition requires that the local behavior of
the outer expansion (\ref{mfpt_b:outex}) as $\x\to\x_i$ must agree
with the far-field behavior as $|\y|\to\infty$ of the inner expansion
(\ref{mfpt_b:innex}), so that
\begin{equation} \label{mfpt_b:mat_1}
  \begin{split}
  \frac{U_0}{\eps} + U_1 + \eps \log\left( \frac{\eps}{2} \right)
  U_2 + &\eps U_3 + \ldots \\
   \sim \frac{U_0}{\eps} & \left( 1 - \frac{C_i}{|\y|} - \frac{\DT_i {\bf \cdot}
       \y}{|\y|^3} \right)  + 
  \log\left( \frac{\eps}{2} \right) V_{1i} + V_{2i} + \ldots \,,
  \end{split}
\end{equation}
where $C_i=C_i(\kappa_i)$ is the reactive capacitance of the
$i$-th patch. Since $|\y|\sim\eps^{-1}|\x-\x_i|$ from
(\ref{app_g:change}) of Appendix \ref{app_g:geod}, we require that
$U_1$ must satisfy (\ref{mfpt_b:Uk}), with the singular behavior
$U_{1} \sim - {U_0 C_{i}/|\x-\x_i|}$ as $\x\to \x_i$ for
$i=1,\ldots,N$, so that
\bsub \label{mfpt_b:U1prob}
\begin{gather}
  \Delta_{\x} U_{1}= 0 \,, \quad \x\in \R^3\backslash \Omega \,; \quad
  \partial_n U_1=0 \,, \quad \x\in \partial\Omega\backslash
  \lbrace{\x_1,\ldots,\x_N\rbrace} \,, \\
  U_1\sim -\frac{U_0 C_{i}}{|\x-\x_i|}\,,  \quad \mbox{as} \quad
  \x\to\x_i \in \partial\Omega \,, \quad i=1,\ldots,N \,,\\
  \lim_{R\to\infty} \int_{\partial\Omega_R} \partial_{r}U_{1}\vert_{r=R} \, ds=-4\pi \,. 
\end{gather}
\esub

From the divergence theorem, the solvability condition for
(\ref{mfpt_b:U1prob}) determines $U_0$ as
\begin{equation} 
  U_0 = -\frac{2}{\overline{C}} \,, \qquad \mbox{where} \qquad
\overline{C}\equiv\sum_{i=1}^{N} C_i(\kappa_i) \,. \label{mfpt_b:U0sol}
\end{equation}
The solution to (\ref{mfpt_b:U1prob}) is represented in terms of the
Green's function of (\ref{mfpt:gs_exact}) as
\begin{equation}
  U_1 = \overline{U}_1 - 2\pi U_0 \sum_{j=1}^{N} C_j G_{s}(\x;\x_j) \,.
  \label{mfpt_b:U1sol}
\end{equation}
As in \cite{Lindsay17}, the unknown constant $\overline{U}_1$ in
(\ref{mfpt_b:U1sol}) has be expanded in terms of additional constants
$\overline{U}_{10}$ and $\overline{U}_{11}$, independent of $\eps$, as
\begin{equation}\label{mfpt_b:swit}
  \overline{U}_1 = \overline{U}_{10} \log\left(\frac{\eps}{2}\right) +
  \overline{U}_{11}\,.
\end{equation}
 The term
$\overline{U}_{10}\log\left({\eps/2}\right)$ in (\ref{mfpt_b:swit}),
which arises from the logarithmic gauge function in
(\ref{mfpt:gs_locm}), is known as a ``switchback term''
(cf.~\cite{Lagerstrom88}), as it effectively inserts a constant term
between ${U_0/\eps}$ and $U_1$ in the outer expansion
(\ref{mfpt_b:outex}). 

To determine $\overline{U}_{10}$ and $\overline{U}_{11}$ we should
proceed to higher order. To do so, we expand $U_1$ in
(\ref{mfpt_b:U1sol}) as $\x\to \x_i$ by using the local behavior
(\ref{mfpt:gs_locm}) of $G_s$ near the $i$-th patch.  The matching
condition (\ref{mfpt_b:mat_1}) becomes
\begin{multline}
  \frac{U_0}{\eps}\left(1 - \frac{C_i}{|\y|}\right) +
  \left(-\frac{U_0 C_i}{2} + \overline{U}_{10}\right)
  \log\left( \frac{\eps}{2} \right) + \frac{U_0 C_i}{2}
 \left(  \frac{y_3 (y_1^2+y_2^2)}{|\y|^3} -\log(y_3+|\y|) \right)\\
 + U_0 \beta_{i} + \overline{U}_{11} + \eps \log\left( \frac{\eps}{2} \right) U_2
 + \eps U_3 + \ldots \\
 \sim \frac{U_0}{\eps}
 \left( 1 - \frac{C_i}{|\y|} - \frac{\DT_i {\bf \cdot}
       \y}{|\y|^3}\right) +\log\left( \frac{\eps}{2} \right)
 V_{1i} + V_{2i} + \ldots \,.
 \label{mfpt_b:mat_2}
\end{multline}
In (\ref{mfpt_b:mat_2}), the constant $\beta_i$ is defined by the
$i$-th component of the matrix-vector product
\begin{equation}\label{mfpt_b:Bi}
  \beta_i =  -2\pi \left({\mathcal G}_s \vc\right)_{i} \,,
\end{equation}
where $\vc\equiv(C_1,\ldots,C_N)^T$, with $C_i=C_i(\kappa_i)$, and
${\mathcal G}_{s}$ is the symmetric Green's matrix,
\bsub \label{mfpt_b:green_mat}
\begin{equation}
    {\mathcal G}_s \equiv \left ( 
\begin{array}{cccc}
 0 & G_{12} & \cdots & G_{1N} \\
 G_{21} & 0 & \cdots   &G_{2N} \\
 \vdots & \vdots  &\ddots  &\vdots\\ 
 G_{N1} &\cdots & G_{N,N-1} & 0
\end{array}
\right ) \,,
\end{equation}
with
\begin{equation}\label{mfpt_b:green_mat_1}
   G_{ij} \equiv G_{s}(\x_i;\x_j) = \frac{1}{2\pi|\x_i-\x_j|} -
  \frac{1}{4\pi}\log\left( 1 + \frac{2}{|\x_i-\x_j|} \right) \,.
\end{equation}
\esub
We remark that since the regular part of the Green's
  function is $-(4\pi)^{-1}\log{2}$ at each patch location $\x_i$, our
  choice of the $\log\left({\eps/2}\right)$ logarithmic gauge function
  in (\ref{mfpt_b:outex}) enforces a zero diagonal in the Green's
  matrix in (\ref{mfpt_b:green_mat}).
  
Upon comparing the ${\mathcal O}(\log\eps)$ terms on both sides of
(\ref{mfpt_b:mat_2}) we conclude that the inner correction $V_{1i}$,
satisfying (\ref{mfpt_b:Vk}) with $k=1$, must
have the far-field behavior $V_{1i}\sim
\overline{U}_{10} - {U_0 C_i/2}$ as $|\y| \to \infty$. As a result,
the solution for $V_{1i}$ is
\begin{equation}
  V_{1i} = \left( \overline{U}_{10} -\frac{U_0 C_i}{2} \right)
  \left(1 - w_{i}\right) \,, \label{mfpt_b:V1sol}
\end{equation}
where $w_{i}=w_{i}(\y;\kappa_i)$ satisfies (\ref{mfpt:wc}). Upon
using the far-field behavior (\ref{mfpt:wc_4}) for $w_{i}$, we obtain that
\begin{equation}
    V_{1i} \sim \left( \overline{U}_{10} - \frac{U_0 C_i}{2}\right) \left(1 - 
      \frac{C_i}{|\y|} +\cdots  \right)\,,
    \quad \mbox{as} \quad |\y| \to \infty\,.
        \label{mfpt_b:V1ff}
\end{equation}

By substituting (\ref{mfpt_b:V1ff}) into the matching condition
(\ref{mfpt_b:mat_2}), we conclude that $U_2$ must
satisfy (\ref{mfpt_b:Uk}) with the singular behavior
$U_2 \sim  -\left(\overline{U}_{10} - \tfrac12 U_0 C_i \right)
{C_i/|\x-\x_i|}$ as $\x\to \x_i$. As a result, we find that $U_2$ satisfies
\bsub \label{mfpt_b:U2prob}
\begin{gather}
  \Delta_{\x} U_{2} = 0 \,, \quad \x\in \R^3\backslash \Omega \,; \quad
  \partial_n U_2=0 \,, \quad \x\in \partial\Omega\backslash
  \lbrace{\x_1,\ldots,\x_N\rbrace} \,, \\
  U_2\sim -\left( \overline{U}_{10} - \frac{U_0 C_i}{2}  \right) \frac{C_i}
  {|\x-\x_i|}\,, \quad \mbox{as} \quad
  \x\to\x_i \in \partial\Omega \,, \quad i=1,\ldots,N \,, \\
 \lim_{R\to\infty} \int_{\partial\Omega_R} \partial_{r}U_{2}\vert_{r=R} \, ds=0  \,. 
\end{gather}
\esub

By using the divergence theorem, we readily find that (\ref{mfpt_b:U2prob}) is
solvable only when $\sum_{j=1}^{N}C_j\left[\overline{U}_{10}-{U_0C_i/2}\right]=0$,
which determines $\overline{U}_{10}$ as
\begin{equation}
  \frac{\overline{U}_{10}}{U_0} = \frac{1}{2\overline{C}} \sum_{j=1}^{N}
  C_j^2   \quad \Rightarrow \quad
  \overline{U}_{10}=-\frac{\sum_{j=1}^{N}C_j^2}{
    \left(\overline{C}\right)^2} \,, \label{mfpt_b:u10}
\end{equation}
where $\overline{C}$ was defined by (\ref{mfpt_b:U0sol}).
With $\overline{U}_{10}$ determined in this way, the solution to
(\ref{mfpt_b:U2prob}) is given in terms of the Green's function in
(\ref{mfpt:gs_exact}) and an additional unknown constant
$\overline{U}_{2}$ as
\begin{equation}
  U_2 = \overline{U}_{2} + 2\pi \sum_{j=1}^{N} C_j \left( \frac{U_0 C_j}{2} -
    \overline{U}_{10} \right)  G_{s}(\x;\x_j)  \,.
  \label{mfpt_b:U2}
\end{equation}

To determine $\overline{U}_{11}$, we must match the ${\mathcal O}(1)$
terms in (\ref{mfpt_b:mat_2}). We obtain that $V_{2i}$ satisfies
(\ref{mfpt_b:Vk}) with $k=2$ subject to the far-field behavior
\begin{equation}
  V_{2i} \sim \beta_i U_0 + \overline{U}_{11} + \frac{U_0 C_i}{2} \left(
     \frac{y_3 (y_1^2 + y_2^2)}{|\y|^3} - \log(y_3 + |\y|)  \right)\,, \quad
   \mbox{as} \quad |\y| \to \infty \,. 
\label{mfpt_b:V2inf}
\end{equation}
Since $V_{0i}=U_0(1-w_{i})$ from (\ref{mfpt_b:v0sol}), we decompose
$V_{2i}$ as 
\begin{equation}\label{mfpt_b:V2decom_1}
  V_{2i} = U_0 \left( \Phi_{2i} +
  \left(\beta_{i} + \frac{\overline{U}_{11}}{U_0}\right)(1-w_{i})\right) \,,
\end{equation}
and obtain from (\ref{mfpt_b:Vk}) and (\ref{mfpt_b:V2inf}) that $\Phi_{2i}$
satisfies
\bsub \label{mfpt_b:Phi2}
\begin{align}
  \Delta_{\y} \Phi_{2i} &=  \left( 2y_{3} w_{i,y_3 y_3} + 2 w_{i,y_3}
                       \right) , \quad
   \y \in \R_{+}^{3} \,, \label{mfpt_b:Phi2_1}\\
   -\partial_{y_3} \Phi_{2i} + \kappa_i \Phi_{2i} &=0 \,, \quad y_3=0 \,,\,
    (y_1,y_2)\in \PT_i\,,  \label{mfpt_b:Phi2_2}\\
    \partial_{y_3} \Phi_{2i} &=0 \,, \quad y_3=0 \,,\, (y_1,y_2)\notin \PT_i
                               \,, \label{mfpt_b:Phi2_3}\\
  \Phi_{2i} &\sim  \frac{C_i}{2} \left(
 \frac{y_3 (y_1^2 + y_2^2)}{|\y|^3} -     \log(y_3 + |\y|) \right)\,, \quad
   \mbox{as} \quad |\y| \to \infty \,. \label{mfpt_b:Phi2_4} 
\end{align}
\esub

From the analysis of the solution to (\ref{mfpt_b:Phi2}) given in
Appendix \ref{app_h:inn2}, we identify the monopole coefficient
$E_i=E_{i}(\kappa_i)$ from the following refined far-field behavior:
\begin{equation}\label{mfpt_b:Phi2_ff}
  \Phi_{2i} - \frac{C_i}{2} \left(
     \frac{y_3 (y_1^2 + y_2^2)}{|\y|^3} - \log(y_3 + |\y|) \right) \sim
  -\frac{E_i}{|\y|}\,, \quad \mbox{as} \quad |\y| \to \infty \,.
\end{equation}
For an arbitrary patch shape $\PT_i$, the determination of
$E_i(\kappa_i)$ is reduced to quadrature in (\ref{app_h:eval}) of
Appendix \ref{app_h:inn2}. We recall that some properties of
$E_i(\kappa_i)$ for circular patches were summarized in Lemma
\ref{lemma:Ej_kappa} of \S \ref{prel:high}, with a highly accurate but
heuristic formula for $E_i=E_i(\kappa_i)$ given in
(\ref{mfpt:E_heur}).

To determine $\overline{U}_{11}$, we will impose a solvability
condition for the problem for the outer correction $U_3$ in
(\ref{mfpt_b:outex}).  To derive the problem for $U_3$, we substitute
(\ref{mfpt_b:Phi2_ff}) into (\ref{mfpt_b:V2decom_1}) and use
$w_{i}\sim {C_i/|\y|}$ as $|\y|\to \infty$. We conclude that $V_{2i}$
satisfies the refined far-field behavior
\begin{equation}\label{mfpt_b:V2ff_refine}
  \begin{split}
    V_{2i} &\sim \beta_i U_0 + \overline{U}_{11} + \frac{U_0 C_i}{2} \left(
      \frac{y_3 (y_1^2 + y_2^2)}{|\y|^3}-  \log(y_3 + |\y|) \right) \\
    &\qquad - \left(  E_i U_0 + \left(\beta_i U_0 + \overline{U}_{11}
      \right)C_i\right) \frac{1}{|\y|} \,, \quad   \mbox{as}
    \quad |\y| \to \infty \,. 
  \end{split}
\end{equation}
The monopole term in (\ref{mfpt_b:V2ff_refine}), given by the
coefficient of ${1/|\y|}$, is one of the two terms that needs to be
accounted for by $U_3$ in the matching condition (\ref{mfpt_b:mat_2}).
The second term is the dipole term in (\ref{mfpt_b:mat_2}), which
arises from (\ref{mfpt:wc_4}). This term is written in terms of outer
variables using (\ref{app_g:change}) of Appendix \ref{app_g:geod}.

In this way, we conclude from (\ref{mfpt_b:Uk}), (\ref{mfpt_b:mat_2})
and (\ref{mfpt_b:V2ff_refine}) that $U_3$ must satisfy 
\bsub
\label{mfpt_b:U3prob}
\begin{align}
&  \Delta_{\x} U_{3} = 0 \,, \quad \x\in \R^3\backslash \Omega \,; \quad
             \partial_n U_3=0 \,, \quad \x\in \partial\Omega\backslash
                      \lbrace{\x_1,\ldots,\x_N\rbrace} \,, \\
&  U_3 \sim -\frac{\left[U_0 E_i +\left(\beta_i U_0 +\overline{U}_{11}\right)C_i
        \right]}{|\x-\x_i|} -U_0 \frac{\DT_i {\bf \cdot}
        {\mathcal Q}_i^T (\x-\x_i)}{|\x-\x_i|^3}\nonumber\\
 & \qquad\qquad\qquad \mbox{as} \quad \x\to\x_i \in \partial\Omega \,, \quad
   i=1,\ldots,N \,,\\
 &  \lim_{R\to\infty} \int_{\partial\Omega_R} \partial_{r}U_{3}\vert_{r=R} \, ds=0  \,. 
\end{align}
\esub Here $\DT_i$ is the dipole vector in (\ref{mfpt:wc_4}), while
the orthogonal matrix ${\mathcal Q}_i$ is defined in
(\ref{app_g:change}) in terms of the basis vectors of the geodesic
coordinate system.

By applying the divergence theorem, we find that (\ref{mfpt_b:U3prob}) is
solvable only when 
\begin{equation}\label{mfpt_b:u11_1}
  \frac{\overline{U}_{11}}{U_0}\overline{C} + \sum_{j=1}^{N}
  \beta_j C_j + \overline{E}=0 \,, \qquad \mbox{where} \qquad
  \overline{E}\equiv \sum_{j=1}^{N}E_j \,.
\end{equation}
We remark that the contribution to the solvability condition from the
dipole term vanishes identically by symmetry owing to the fact that
$\DT_i$ has the form $\DT_i=(p_{1i},p_{2i},0)^T$.  Upon recalling
(\ref{mfpt_b:Bi}) for $\beta_i$, we solve (\ref{mfpt_b:u11_1}) for
$\overline{U}_{11}$ to obtain
\begin{equation}\label{mfpt_b:u11}
  \frac{\overline{U}_{11}}{U_0} = \frac{2\pi}{\overline{C}} \vc^T {\mathcal G}_s
  \vc - \frac{\overline{E}}{\overline{C}} \quad \Rightarrow \quad
  \overline{U}_{11}=\frac{2}{\left(\overline{C}\right)^2} \left(
    \overline{E} - 2\pi\vc^T {\mathcal G}_s
  \vc\right) \,.
\end{equation}

Finally, to determine an expression for the capacitance $C_{\rm T}$ in
(\ref{bp:ssp_3}), we must take the limit as $|\x|\to\infty$ of our
outer asymptotic expansion (\ref{mfpt_b:outex}) and compare it with the
required limiting behavior (\ref{bp:ssp_3}). By comparing the
${\mathcal O}(1)$ terms in the resulting expression we obtain a
three-term asymptotic expansion for $C_{\rm T}$. We summarize our main
result for $C_{\rm T}$ as follows.

\begin{prop}\label{mfpt_b:main_res} 
  For $\eps \to 0$, the capacitance $C_{\rm T}$ for the dimensionless
  problem (\ref{mfpt:ssp}) outside the unit sphere in the presence of
  $N$ well-separated partially reactive Robin patches, centered at
  $\x_i$ for $i=1,\ldots,N$ and with local reactivities $\kappa_i$
  for $i=1,\ldots,N$, has the three-term asymptotics
  \bsub \label{mfpt_b:main_res_1}
\begin{equation}
  \frac{1}{C_{\rm T}} \sim \frac{|U_0|}{\eps} \left[ 1 +
    \eps \frac{\overline{U}_{10}}{U_0} \log\left(\frac{\eps}{2}\right)
    + \eps \frac{\overline{U}_{11}}{U_0} + {\mathcal O}\left(\eps^2
    \left[\log\left({\eps/2}\right)\right]^2\right)\right]\,,
\end{equation}
where
\begin{equation}\label{mfpt_b:main_res_c}
  |U_0| = \frac{2}{\overline{C}} \,, \qquad
 \frac{\overline{U}_{10}}{U_0} = \frac{1}{2\overline{C}} \sum_{i=1}^{N}
 C_i^2    \,, \qquad  \frac{\overline{U}_{11}}{U_0} =
 \frac{2\pi}{\overline{C}} \vc^T {\mathcal G}_s
  \vc - \frac{\overline{E}}{\overline{C}}\,.
\end{equation}
\esub Here $\overline{C}=\sum_{i=1}^{N}C_i$ and
$\overline{E}=\sum_{i=1}^{N}E_i$.  Properties of the reactive
capacitance $C_i=C_i(\kappa_i)$ and the monopole coefficient
$E_i=E_i(\kappa_i)$, as defined by the far-field behaviors in
(\ref{mfpt:wc_4}) and (\ref{mfpt:inn2_h4}), respectively, were
summarized in Lemmas \ref{lemma:Cj_kappa} and
\ref{lemma:Ej_kappa}. The Green's matrix ${\mathcal G}_s$ in
(\ref{mfpt_b:main_res_c}) is given by (\ref{mfpt_b:green_mat}). In
terms of the dimensional capacitance ${\mathcal C}_{\rm T}$, defined in
(\ref{mfpt:ssp0_c}), we use (\ref{intro:scalings}) to obtain for a
sphere of radius $R$ and with a collection of partially reactive
patches of maximum diameter $L$ that
\begin{equation}\label{mfpt_b:dimen}
  {\mathcal C}_{\rm T} = R C_{\rm T} \,.
\end{equation}
Here in calculating $C_{\rm T}$ from (\ref{mfpt_b:main_res_1}) we use
$C_i=C_i\left({L\K_i/D}\right)$ and $E_i=E_i\left({L\K_i/D}\right)$ and
 we evaluate the Green's matrix ${\mathcal G}_s$ at $\x_i={\X_i/R}$.
\end{prop}

We remark that for the special case of $N$ locally circular patches
with perfect reactivity ($\kappa_i=\infty$), for which
$C_i={2a_i/\pi}$ and
$E_i=-2 a_i^2{\left(\log{a_i} + \log{4} - {3/2} \right)/\pi^2}$, our
three-term asymptotic result (\ref{mfpt_b:main_res_1}) for
${1/C_{\rm T}}$ agrees with that in Principal Result 3.1 of
\cite{Lindsay17}. Our main result (\ref{mfpt_b:main_res_1}) extends
the previous result of \cite{Lindsay17} to arbitrary-shaped patches
and to finite reactivity $\kappa_i>0$. We emphasize that, as
  expected intuitively, our result is identical in form to that
  derived in \cite{Lindsay17} provided that we simply replace the two
  monopole coefficients in Eq. (3.37) of \cite{Lindsay17} with
  $C_i(\kappa_i)$ and $E_{i}(\kappa_i)$.  Moreover, we remark that
our three-term asymptotic result remains well-ordered in $\eps$ over
the full range of patch reactivities. As a result, we can explore the
limits $\kappa_i\ll 1$ and $\kappa_i\gg 1$ directly from
(\ref{mfpt_b:main_res_1}).

There are a few additional key features of our main result for $C_{\rm T}$.
Firstly, since $\overline{U}_{11}$ depends on the Green's matrix
${\mathcal G}_s$, we conclude that only the third-order term in the
asymptotic expansion of $C_{\rm T}$ depends on the spatial configuration
$\lbrace{\x_1,\ldots,\x_N\rbrace}$ of the centers of the reactive
patches on the surface of the unit sphere. It is this term that
incorporates spatial effects from the diffusive interactions between
patches. Secondly, since the leading-order term in
(\ref{mfpt_b:main_res_1}) does not depend on the details of the
Green's function, it is valid for an arbitrary bounded domain
${\mathcal B}$ with smooth boundary. By retaining only the
leading-order terms in our three-term main result
(\ref{mfpt_b:main_res_1}), we obtain
\begin{equation}
  {\mathcal C}_{\rm T} \approx \sum\limits_{i=1}^N
  \frac{L}{2} C_i(L \K_i/D) \,.
\end{equation}
This shows that the capacitance of the structured target is, to
leading order, the sum of reactive capacitances $LC_i/2$ of individual
patches, in analogy to a parallel connection of capacitors in
electrostatics.  The factor $1/2$ accounts for the fact that only one
side of the patch located on the boundary is accessible to Brownian
particles.  As in the Berg-Purcell result, this leading-order term
does not account for diffusive interactions between the patches (which
appears only in the third-order term), nor for the curvature of the
boundary.
Thirdly, to numerically evaluate the asymptotic result
(\ref{mfpt_b:main_res_1}), we need only compute $C_i(\kappa_i)$ from
the Steklov eigenfunction expansion (\ref{eq:Cmu_def0}) and
$E_i(\kappa_i)$ from the quadrature in (\ref{app_h:eval}) of Appendix
\ref{app_h:inn2}. When the Robin patches are all locally circular, by
using the heuristic approximations given in (\ref{mfpt:sigmoidal_2})
and (\ref{mfpt:E_heur}) for $C_i$ and $E_i$, respectively, we obtain
an explicit expression for $C_{\rm T}$ valid for all $\kappa_i>0$.

\section{The Effective Capacitance and Effective Reactivity}\label{sec:homog}

We now derive a new scaling law for the effective capacitance $\ceff$
and the effective reactivity $\keff$ for (\ref{mfpt:ssp}) applicable
to a large number (i.e. $N\gg 1$) of uniformly distributed identical
circular patches of a common radius $\eps$ and reactivity
$\kappa$. Our result applies to the low patch area fraction
  limit where $f\equiv {N\pi \eps^2/|\partial\Omega|}$, with
  $|\partial\Omega|=4\pi$, satisfies
  ${\mathcal O}(\eps^2)\ll f \ll 1$. A similar result was given in
\cite{Lindsay17} for perfectly reactive patches, and in
\cite{Cheviakov13} for the mean first-passage time (MFPT) narrow
capture problem within a sphere with small surface patches of finite
reactivity.

For the case of $N$ identical circular patches, we set $a_i=1$ for
$i=1,\ldots,N$ in (\ref{mfpt_b:main_res_1}) to obtain that our
three-term expansion becomes
\begin{equation}\label{homo:c0_orig_1}
  \frac{1}{C_{\rm T}} \sim \frac{2}{N C \eps} \left[ 1 + \frac{\eps C}{2}
    \log\left(\frac{\eps}{2}\right) + \eps C\left( \frac{2}{N}
      {\mathcal H}(\x_1,\ldots,\x_N) - \frac{E}{C^2} \right)\right]\,.
\end{equation}
Here $C=C(\kappa)$ is the common reactive capacitance of each patch,
defined by the far-field behavior of the solution $w = w(\y;\kappa)$
to (\ref{mfpt:wc}), where $\Gamma$ is the unit disk and where we
omit the subscript $i$. In addition, the monopole coefficient
$E=E(\kappa)$ in (\ref{homo:c0_orig_1}) is given by
(\ref{mfpt:Ej_all}) in which we set $a_i=1$ and omit the subscript $i$.
In (\ref{homo:c0_orig_1}), the discrete energy ${\mathcal H}$,
representing inter-patch interactions from the Green's
matrix, is
\begin{equation}\label{homo:H}
    {\mathcal H}(\x_1,\ldots,\x_N)=\sum_{i=1}^{N}\sum_{j=i+1}^{N}
    \left( \frac{1}{|\x_i-\x_j|} - \frac{1}{2}\log\left(1 +
        \frac{2}{|\x_i-\x_j|}\right) \right) \,.
\end{equation}

As discussed in \cite{Lindsay17} (see also analogous results in
Appendix of \cite{Cheviakov13} and \cite{Cheviakov10} for the MFPT
narrow capture problem), for a large collection of uniformly
distributed patches on the surface of the sphere with centers at
$\x_i$ for $i=1,\ldots,N$, a mean-field approximation for the
discrete energy ${\mathcal H}$ has the scaling law
\begin{equation} \label{homo:hasy}
  {\mathcal H} = \frac{N^2}{4} - d_1 N^{3/2} + 
  \frac{1}{8} N\log{N} + d_2 N + d_3 N^{1/2} + {\mathcal O}(\log{N}) \,,
  \quad \mbox{for} \quad N\gg 1\,,
\end{equation}
for some coefficients $d_1$, $d_2$, and $d_3$. Such a
  scaling law, should also correspond to the global minimal
  energy for ${\mathcal H}$.  As discussed in \S 4 of
  \cite{Lindsay17}, for points on either a sphere or a flat surface,
  the choice $d_1\approx 0.55230$ is the best fit to the global
  minimum of ${\mathcal H}$ (see also \cite{Rakhamanov94}), and for a sphere
  corresponds roughly to a uniform arrangement of points. The further choices
  $d_2 = 1/8$ and $d_3 = 1/4$ result from a formal mean-field
  approximation for uniformly distributed points (see \S 4 of
  \cite{Lindsay17} for a further discussion).

Defining $\zeta=\zeta(\kappa)\equiv e^{-2E/C^2}e^{4d_2}$, the
effective capacitance $\ceff$ is obtained by substituting
(\ref{homo:hasy}) into (\ref{homo:c0_orig_1}), which yields
\begin{equation}\label{homo:c0_2}
  \frac{1}{\ceff}\sim \frac{2}{NC \eps} + 1 +\frac{1}{N}
  \log\left( \frac{\eps\sqrt{N}}{2}\zeta \right) - \frac{4d_1}{N^{1/2}}
  + \frac{4d_3}{N^{3/2}} +
  {\mathcal O}\left(\frac{\log{N}}{N^2}\right)\,.
\end{equation}

Next, we write (\ref{homo:c0_2}) in terms of the patch area fraction
$f\equiv {N\eps^2/4}$. Since for the well-ordering of our asymptotic
expansion the first term in (\ref{homo:c0_2}) must be the dominant
term, we require that $N\ll {\mathcal O}(\eps^{-1})$. This enforces
that our result will apply to the small patch area fraction limit
$f\ll 1$. Upon eliminating $N$ in (\ref{homo:c0_2}) using
$N={4f/\eps^2}$, and assuming that ${\mathcal O}(\eps^2)
\ll f \ll 1$, we obtain our main result for the dimensionless
effective capacitance:
\bsub \label{homo:ceff}
\begin{equation}\label{homo:ceff_1}
  \frac{1}{\ceff} \sim 1 + \frac{\eps}{2Cf} \left[ 1 - 4d_1
    C \sqrt{f} + \frac{\eps C}{2}\log\left(\zeta \sqrt{f}\right)
    + \frac{d_3 C \eps^2}{\sqrt{f}} +
    {\mathcal O}\left(\frac{\eps^3}{f}
        \log\left({\sqrt{f}/\eps}\right) \right) \right]\,,
\end{equation}
where, we have
\begin{equation}\label{homo:ceff_2}
  f= \frac{\eps^2 N}{4}\ll 1 \,, \quad
  \zeta(\kappa) =e^{-2E/C^2}e^{4d_2} \,, \quad d_1\approx 0.5523\,, \quad
   d_2=\frac{1}{8} \,, \quad d_3=\frac{1}{4} \,.
\end{equation}
\esub The dependence of $\ceff$ on the reactivity $\kappa$ is
inherited through $C(\kappa)$ and $\zeta(\kappa)$.

To determine the effective reactivity, we use the dimensionless form of
(\ref{eq:kappa_J}) given by ${1/\keff}=-1+ {1/\ceff}$.
In this way, by using (\ref{homo:ceff_1}) for $\ceff$ in
this result, we obtain a scaling law for the dimensionless
effective reactivity given by
\begin{equation}\label{homo:keff}
\keff  \sim \frac{2C f}{\eps} \left[ 1 - 4d_1
    C \sqrt{f} + \frac{\eps C}{2}\log\left(\zeta(\kappa) \sqrt{f}\right)
    + \frac{d_3 C \eps^2}{\sqrt{f}} +
        {\mathcal O}\left(\frac{\eps^3}{f}
        \log\left({\sqrt{f}/\eps}\right) \right) \right]^{-1}\,.
\end{equation}
Setting $d_3={1/4}$, and neglecting the error term in
  (\ref{homo:keff}), we can also rewrite this expression in an
  alternative form as
\begin{equation}  \label{eq:keff}
  \keff \sim \frac{N \, \eps C(\kappa)}{2} \left[1 - \sqrt{N} \, \eps
  C(\kappa)\left(2d_1 - \frac{1}{2N} 
 -\frac{1}{2\sqrt{N}}
\log\bigl(\zeta(\kappa)\eps\sqrt{N}/2 \bigr)\right)\right]^{-1} \,.
\end{equation}
Even though (\ref{homo:keff}) or (\ref{eq:keff}) can be formally
re-expanded by using the binomial approximation, the resulting
expressions turn out to be less accurate than (\ref{homo:keff}) or
(\ref{eq:keff}) when $\eps$ is not too small, as witnessed by
comparison with Monte Carlo simulations.  For this reason, we avoid
using such re-expansions and retain (\ref{homo:keff}) or
(\ref{eq:keff}) in the following numerical illustrations.

Our main results for $\ceff$ and $\keff$ are determined in terms of
both $C(\kappa)$ and $E(\kappa)$. For $\kappa\ll 1$, we set $a_i=1$ in
(\ref{mfpt:cj_small_b}) and (\ref{esmall}) to obtain
\begin{equation}\label{homo:low}
  C(\kappa)\sim \frac{\kappa}{2} -\frac{4}{3\pi}\kappa^2 + 0.3651 \,
  \kappa^3 \,, \qquad
  E(\kappa)\sim \frac{\kappa^2}{32} \quad \mbox{as} \quad \kappa\to 0\,.
\end{equation}
To the leading order in $\kappa$, this yields ${E/C^2}\sim {1/8}$ so
that 
$\zeta\sim e^{1/4} \approx 1.2840$ as $\kappa\to 0$. In
contrast, in the limit $\kappa\to\infty$ of large reactivity, we have
from setting $a_i=1$ in (\ref{mfpt:cj_large_b}) and
(\ref{mfpt:EJ_asy_1}) that
\begin{equation}\label{homo:high}
    C\sim {2/\pi} \,, \qquad E\sim (3 - 4\log 2)/\pi^2 \quad \mbox{as} \quad \kappa \to \infty \,.
\end{equation}
This gives ${E/C^2}\sim (3-4\log 2)/4$, which yields
that $\zeta\sim 4/e \approx 1.4715$ as $\kappa\to \infty$ in
(\ref{homo:keff}). In this way, we recover the scaling law derived in
\S 4 of \cite{Lindsay17} for the effective reactivity for the case of
perfectly reactive patches.
  
In summary, the results (\ref{homo:ceff}), (\ref{homo:keff}) and
(\ref{eq:keff}) determine the dimensionless effective capacitance and
effective reactivity over the full range $\kappa>0$ of reactivity of a
large collection of identical uniformly distributed circular patches
of a common radius $\eps$. By setting $a_i=1$ in the heuristic
approximations for $C$ and $E$ in (\ref{mfpt:sigmoidal_2}) and
(\ref{mfpt:E_heur}), we obtain explicit analytical results for these
two effective parameters for all $\kappa>0$.

Next, we deduce the corresponding homogenized result for the
dimensional problem (\ref{mfpt:ssp0}) in a sphere of radius $R$ and
with circular patches of a common radius $L\ll R$ and common
dimensional reactivity $\K_i=\K$ for $i=1,\ldots,N$. Upon recalling
the scaling relations (\ref{intro:scalings}), we have from
(\ref{homo:keff}) that
\begin{equation}\label{homo:dimen}
    \Keff = \frac{D}{L} \keff = \frac{D}{\eps R} \keff \,,
\end{equation}
where $C=C\left({L\K/D}\right)$ and $E=E\left({L\K/D}\right)$.
The dimension of $\Keff$ is ${\mbox{length}/\mbox{time}}$.

Finally, in order to compare with the heuristic approximation
(\ref{eq:kapheur_2}) of \cite{Berez04}, we consider our leading-order
theory from (\ref{homo:keff}) where we use $\keff\sim {2C/f}$ together
with $L=\eps R$ for the patch radius and (\ref{mfpt:sigmoidal_2}) to
approximate $C(\kappa)$. This yields that
\begin{equation}\label{homo:dimen_lead}
  \Keff = \frac{1}{\eps} \left( \frac{2 fD}{\eps R}\right)
  {\mathcal C}\left(\frac{\K \eps}{D}\right) \,, \qquad \mbox{where} \qquad
  {\mathcal C}(\mu)\equiv \frac{{2\mu/\pi}}{\mu + {4/\pi}} \,,
\end{equation}
with the extra pre-factor of ${1/\eps}$ resulting from the fact that
we homogenized $\eps \partial_n U + \kappa_iU=0$ in (\ref{bp:ssp_2}) rather
than $\partial_n U + \kappa_iU=0$. We conclude that our leading-order
theory exactly reproduces the heuristic formula (\ref{eq:kapheur_2}) of
\cite{Berez04}, which was first derived theoretically in \cite{Plunkett24}
as was discussed in \S \ref{all:intro}.

\section{Monte Carlo simulations}\label{sec:numer}

To assess the accuracy and the range of validity of our asymptotic
results, we resort to a numerical calculation of the capacitance
${\mathcal C}_{\rm T}$.  For this purpose, one can generate random
trajectories of diffusing particles and determine the flux $J$, from
which the effective reactivity $\K_{\rm eff}$ and the capacitance
${\mathcal C}_{\rm T}$ follow from (\ref{eq:kappa_J}).  In the steady-state
regime, any trajectory terminates either by reaction on a patch, or by
escape to infinity (due to the transient character of diffusion in
three dimensions), and the fraction of reacted particles is
proportional to the flux (see below).  Different numerical schemes
were proposed to undertake such Monte Carlo simulations, including
variable-jump techniques based on the {\em walk-on-spheres} algorithm
\cite{Zhou17,Schumm23,Ye25}.  However, for all of these approaches,
modeling reflections on an inert boundary is the most difficult and
time-consuming part of the computation.  We now propose an alternative,
highly efficient Monte Carlo algorithm, which relies on the rotational
symmetry of the spherical domain.

\subsection{Efficient algorithm for a sphere} 

Let us introduce a spherical surface of radius $R+\ell$,
$\pa^\ell = \{\X\in\R^3 \, \vert \,\, |\X| = R+\ell\}$, and
define the sequence of stopping times for a random trajectory $\X_t$
of a diffusing particle as
\begin{equation}
\tau_{2j+1} = \inf\{ t > \tau_{2j} ~:~ \X_t \in \pa\}\,,
\qquad \tau_{2j} = \inf\{ t > \tau_{2j-1} ~:~ \X_t \in \pa^\ell\}\,,
\end{equation} 
with $j=1,2,\cdots$, and $\tau_0 \equiv 0$.  In other words, $\tau_1$
is the first-passage time to the sphere $\pa$, $\tau_2$ is the first
instance when $\X_t$ crosses $\pa^\ell$ after $\tau_1$,
$\tau_3$ is the first instance of the return to the sphere after
$\tau_2$, etc.  The stopping times $\tau_j$ partition the trajectory
into consecutive {\it independent} paths:
$\X_0 \rightsquigarrow \X_{\tau_1}$,
$\X_{\tau_1} \rightsquigarrow \X_{\tau_2}$, etc.  If the particle
escapes to infinity (and thus never returns to $\pa$) or reacts on
$\pa$ during the $j$-th path, all stopping times $\tau_i$ with
$i\geq j$ are set to $+\infty$.

(i) The first path consists in the first arrival onto the sphere from
a given starting point $\X_0$.  Instead of simulating the random
trajectory over $t \in (0,\tau_1)$, one can generate the first arrival
point $\X_{\tau_1} = \hat{\X}_1$ on the sphere according to the
harmonic measure density, which is given explicitly for the exterior
of a sphere of radius $R$ by
\begin{equation}  \label{eq:HM}
\omega_{\X_0}(\X_1) = \frac{|\X_0|^2 - R^2}{4\pi R\, |\X_1 - \X_0|^3} \,.
\end{equation}
Let us introduce the local spherical coordinates with the north pole
oriented along $\X_0$ such that $\X_0^{\prime} = (r_0,0,0)$, while
$\X^{\prime}_1 =(R,\theta,\phi)$ is determined by the angles $\theta$
and $\phi$, where the prime indicate the local coordinates.  By
symmetry, the azimuthal angle $\phi$ is uniformly distributed over
$(0,2\pi)$.  In turn, the probability density of the polar angle
$\theta \in (0,\pi)$ follows from (\ref{eq:HM}) as
\begin{equation}  \label{eq:HM_theta}
  \omega_{\rho}(\theta) = \frac{\rho^2 - 1}{2 (1 - 2\rho \cos\theta +
    \rho^2)^{3/2}} \sin\theta  \,,
\end{equation}
where $\rho = {r_0/R} > 1$. Here we multiplied (\ref{eq:HM}) by $2\pi
\sin\theta$ to account for the Jacobian and for the already generated
azimuthal angle.  The integral of this density yields the cumulative
distribution function
\begin{equation}\label{mc:fp}
  F_{\rho}(\theta) = \int\limits_0^\theta  \omega_{\rho}(\theta^{\prime})\,
  d\theta^{\prime}
  = \frac{\rho+1}{2\rho} \biggl(1 -
  \frac{\rho-1}{\sqrt{1 - 2\rho \cos\theta + \rho^2}}\biggr)\,.
\end{equation} 
One easily checks that $F_{\rho}(0) = 0$ and $F_{\rho}(\pi) = 1/\rho <
1$.  As expected, the probability densities in (\ref{eq:HM}) and
(\ref{eq:HM_theta}) are not normalized to unity owing to the transient
character of diffusion in three dimensions.  In fact, the particle can
escape to infinity with probability $1 - 1/\rho = 1 - R/r_0$.  To
account for this possibility, one generates a random variable $\eta$
with uniform distribution on $(0,1)$.  If $\eta > 1/\rho$, the
simulation stops due to escape to infinity.  Otherwise, one sets $\eta
= F_{\rho}(\theta)$ and inverts the relation (\ref{mc:fp}) to generate the
angle $\theta$ of the arrival point as
\begin{equation}
\theta = \arccos\biggl(1 - \frac{(\rho-1)^2}{2\rho} \biggl(\frac{1}{(1 - \frac{2\rho}{1+\rho} \eta)^2} - 1\biggr)\biggr).
\end{equation}  
In this way, one generates the random point $\hat{\X}^{\prime}_1 =
(R,\theta,\phi)$ on the sphere in the local spherical coordinates.  To
complete this step, it remains to transform these local coordinates
into the global coordinates to get $\X_{\tau_1} = \hat{\X}_1$.  We
emphasize that in this step a long simulation of the random path $\X_0
\rightsquigarrow \X_{\tau_1}$ is replaced by the generation of two
random variables $\phi$ and $\theta$ with explicitly known
distributions.

(ii) The next path is diffusion near a partially reactive sphere until
hitting the surface $\pa^\ell$ (or reacting on $\pa$).
This is the most time-consuming and delicate step for common Monte
Carlo techniques, which use a very small timestep to approximate the
random path $\X_{\tau_1} \rightsquigarrow \X_{\tau_2}$ by a sequence
of small random jumps, with eventual reflections/reactions on the
boundary.  Following the concept of the {\em escape-from-a-layer}
approach \cite{Ye25}, we aim at replacing such a detailed and
time-consuming simulation by a single escape event.  Let us first
consider the case when the reactivity is homogeneous. The extension to
the heterogeneous case is discussed below.  With a homogeneous
reactivity, the splitting probability $\Pi(r)$ of the escape event can
be easily found by solving the Laplace equation in a spherical layer
$\Omega^\ell = \{\X\in \R^3 \,\vert \,\, R < |\X| <
  R+\ell\}$, with Dirichlet boundary condition at
$r = R+\ell$ and Robin boundary condition at $r = R$,
yielding
\begin{equation}
  \Pi(r) = 1 - \frac{\K R/D}{1 + \frac{\K \ell/D}{1+\ell/R}}
    \biggl(\frac{R}{r} -
  \frac{R}{R+\ell}\biggr)  \qquad \mbox{for} \quad R \leq r \leq R+\ell\,.
\end{equation}
For the starting point on the reactive sphere, one has
\begin{equation} \label{eq:Pi_escape}
\Pi(R) = \biggl(1 + \frac{\ell\K/D}{1+\ell/R}\biggr)^{-1} \,.
\end{equation}
In other words, one can generate a uniform random variable $\eta \in
(0,1)$ to choose between two possible outcomes: (a) if $\eta >
\Pi(R)$, the particle reacts on $\pa$ before escaping to $\pa^\ell$,
and the simulation stops; or (b) if $\eta < \Pi(R)$, the particle hits
$\pa^\ell$ before reacting on $\pa$.  To complete this step, one needs
to generate the random escape position $\X_{\tau_2} = \hat{\X}_2$ on
$\pa^\ell$.  Even though the conditional distribution of $\hat{\X}_2$
can be derived and implemented, we resort to a much simpler
approximation, which consists in replacing the random variable
$\hat{\X}_2$ by its mean value, which by symmetry is simply
$\hat{\X}_1 (R+\ell)/R$ (i.e., the previous hitting point $\hat{\X}_1$
on the sphere $\pa$ is lifted to $\pa^\ell$ by rescaling).  In fact,
it is unlikely for a particle to diffuse far away from the point
$\hat{\X}_1$ in a thin reactive layer $\Omega^\ell$ so that the
distribution of $\hat{\X}_2$ is peaked around $\hat{\X}_1
(R+\ell)/\ell$ when $\ell/R \ll 1$.  We conclude that the most
time-consuming simulation of the path $\X_{\tau_1} \rightsquigarrow
\X_{\tau_2}$ is replaced by the generation of the uniform random
variable $\eta$.

(iii) From the point $\hat{\X}_2$, the particle resumes its diffusion
until its arrival onto $\pa$, or escape to infinity.  The simulation
of the next path $\X_{\tau_2} \rightsquigarrow \X_{\tau_3}$ can thus
be substituted by the generation of the arrival point $\X_{\tau_3} =
\hat{\X}_3$ according to the harmonic measure density, as described in
step (i).  After the arrival onto the boundary $\pa$, the particle
needs to escape from a thin layer of width $\ell$, as described in
step (ii), and so on.  As a consequence, the detailed simulation of
the random trajectory $\X_t$ is replaced by a succession of steps (i)
and (ii), that sample a sequence of random points $\hat{\X}_1,
\hat{\X}_2, \hat{\X}_3, \cdots$ of that trajectory.  The simulation of
the trajectory is stopped when either the particle escapes to
infinity, corresponding to a step with some odd index $2j-1$, or
reacts on the sphere during a step with some even index $2j$.
Repeating such a simulation $M$ times, one can approximate the
probability of reaction on the sphere as $P_{\rm react}(\X_0) \approx
M_{\rm react}/M$, where $M_{\rm react}$ is the number of simulated
trajectories that terminated by reaction.  Finally, the latter
probability determines the flux of particles started from $\X_0$ as
$J(\X_0) = J_{\rm Smol} P_{\rm react}(\X_0)$.

Since the steady-state flux $J$ from (\ref{eq:J_kappa}) corresponds to
a source at infinity, in the sense that ${\mathcal U}(\X)\to {\mathcal
U}_{\rm inf}$ as $|\X|\to\infty$, one would formally need to start
simulations with $\X_0$ at infinity.  However, the rotational symmetry
of the sphere implies that the first arrival point $\hat{\X}_1$ is
distributed uniformly on the sphere so that one can simply replace the
very first step (the path from $\X_0$ to $\hat{\X}_1$) by the
generation of the uniform point $\hat{\X}_1$ on $\pa$.  Keeping the
remaining steps unchanged, one can thus determine the probability of
reaction $P_{\rm react}(\circ)$ with the uniform starting point, from
which $J = J_{\rm Smol} P_{\rm react}(\circ)$.  Substituting this
expression into (\ref{eq:kappa_J}), one can directly access the
effective reactivity:
\begin{equation}
\frac{R}{D} \K_{\rm eff} = \biggl(\frac{1}{P_{\rm react}(\circ)} - 1\biggr)^{-1} \,.
\end{equation}

\begin{figure}
\begin{center}
\includegraphics[width=88mm]{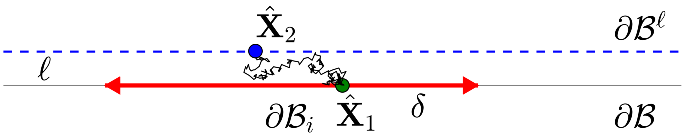} 
\end{center}
\caption{ Schematic two-dimensional illustration of a thin layer of
  width $\ell$ near the boundary $\pa$ with a reactive
  patch $\pa_i$ (thick red interval).  A random path
  $\hat{\X}_1 = \X_{\tau_1} \rightsquigarrow \X_{\tau_2} = \hat{\X}_2$
  from the first arrival point $\hat{\X}_1$ to the escape point
  $\hat{\X}_2$ on the surface $\pa^\ell$ is shown.  Two
  arrows indicate the boundary $\pae_{i}$ of the patch $\pa_i$ (in
  three dimensions, $\pae_{i}$ is a curve but here it is reduced to
  two endpoints of the shown interval).}
\label{fig:scheme1}
\end{figure}

In the present form, the proposed algorithm is formally limited to the
case of homogeneous reactivity, for which the solution is already
known.  In the heterogeneous case, the reactivity is a piecewise
constant function on the spherical boundary, which can be written in
the form
\begin{equation}
  \K(\X) = \sum\limits_{i=1}^N \K_i \, {\mathcal I}_{\pa_i}(\X)  \qquad
  \mbox{for} \quad \X\in \pa\,,
\end{equation}
where ${\mathcal I}_{\pa_i}(\X)$ is the indicator function of the
$i$-th patch $\pa_i$ with reactivity $\K_i$.  As a consequence, the
escape probability from a thin layer depends on the spatial
arrangement of patches and on their reactivities, and
(\ref{eq:Pi_escape}) is in general not applicable.  However, if the
layer width $\ell$ is much smaller than the patch sizes,
one can still employ (\ref{eq:Pi_escape}) by setting $\K$ to be the
reactivity at the arrival point $\hat{\X}_1\in \pa$ (or at
$\hat{\X}_{2j-1}\in \pa$ for other paths).  For instance, let us
assume that $\hat{\X}_1 \in \pa_i$, and we denote
$\delta \equiv |\hat{\X}_1 - \pae_{i}|$ to be the distance between
$\hat{\X}_1$ and the boundary $\pae_{i}$ of the patch $\pa_{i}$ (see
Fig. \ref{fig:scheme1}).  If $\delta \gg \ell$, the
probability of diffusing inside the reactive thin layer
${\mathcal B}_\ell$ from $\hat{\X}_1$ to points farther than the
distance $\delta$ is exponentially small, of the order
${\mathcal O}(e^{-\delta/\ell})$.  In other words, since it
is unlikely that the particle may be exposed to a reactivity other
than $\K_i$, the formula (\ref{eq:Pi_escape}) can be used, and the
algorithm described above applies. Similarly, if $\hat{\X}_1$ does not
belong to any reactive patch, the particle diffuses near an inert
reflecting boundary with reactivity $\K = 0$, and the escape
probability is equal to $1$.  The only difficult case is the situation
when the arrival point $\hat{\X}_1$ occurs in the very close vicinity
of the boundary $\pae_{i}$ of a reactive patch. In this case,
(\ref{eq:Pi_escape}) is not valid because the particle may encounter
the boundary on either side of $\pae_{n}$ that would alter $\Pi(R)$.
However, since the relative surface area of the regions near the
boundary of the patches, estimated as
$2\ell (|\pae_{1}| + \ldots + |\pae_{N}|)/(4\pi R^2) \propto \ell\eps
N/R^2$, is small, the contribution of inaccurate estimations of the
escape probability via (\ref{eq:Pi_escape}) is expected to be small as
well.  This is actually confirmed by our numerical simulations, as
discussed below.
In summary, we still employ (\ref{eq:Pi_escape}) for the heterogeneous
case by setting $\K = \K(\hat{\X}_{2j-1})$.

\subsection{Validation}\label{sec:mc:validate}

For most of our computations below and in \S \ref{sec:comp} we set
$\ell = 10^{-2}$ and used $M = 10^5$ realizations.  We
verified that a further decrease of $\ell$ did not affect
the numerical results.  Although an increase of $M$ would result in
smaller statistical errors, given that the relative error is of the
order of $1/\sqrt{M}$, our choice of $M=10^{5}$ was sufficient to
provide good accuracy.

Our Monte Carlo algorithm has been validated in several settings, as
described in (i)--(iii) below, such as when the solution is either
known analytically or could be obtained numerically by other methods.

(i) When the reactivity is homogeneous, the probability of reaction
$P_{\rm react}(\X_0)$ is known explicitly as $P_{\rm react}(\X_0) =
\tfrac{R}{|\X_0|} \, \frac{1}{1+ D/(\K R)}$.  This relation was
used to validate the algorithm (results not shown).  In addition, one
can quantify the distribution of the reaction events on the sphere via
the spread harmonic measure \cite{Grebenkov06a,Grebenkov15}.  In the
spherical coordinates oriented such that the starting point is $\X_0 =
(r_0, 0, 0)$, the random position $\X_r = (R, \theta,\phi)$ of the
reaction event is characterized by the uniformly distributed azimuthal
angle $\phi$ and the polar angle $\theta$ obeying the following
probability density in terms of the Legendre polynomials $P_{n}(z)$:
\begin{equation}  \label{eq:SHM_theta}
  \omega_{r_0,\K}(\theta) = \sin\theta \sum\limits_{n=0}^\infty
  P_n(\cos\theta) (R/r_0)^{n+1} \, \frac{n+1/2}{1 + (n+1)D/(\K R)} \,.
\end{equation}
In the limit $\K = \infty$, this spread harmonic measure density
reduces to (\ref{eq:HM_theta}).  The comparison of
(\ref{eq:SHM_theta}) to an empirical probability density obtained from
Monte Carlo simulations served for validation purposes (results not
shown).

\begin{figure}
\begin{center}
\includegraphics[width=88mm]{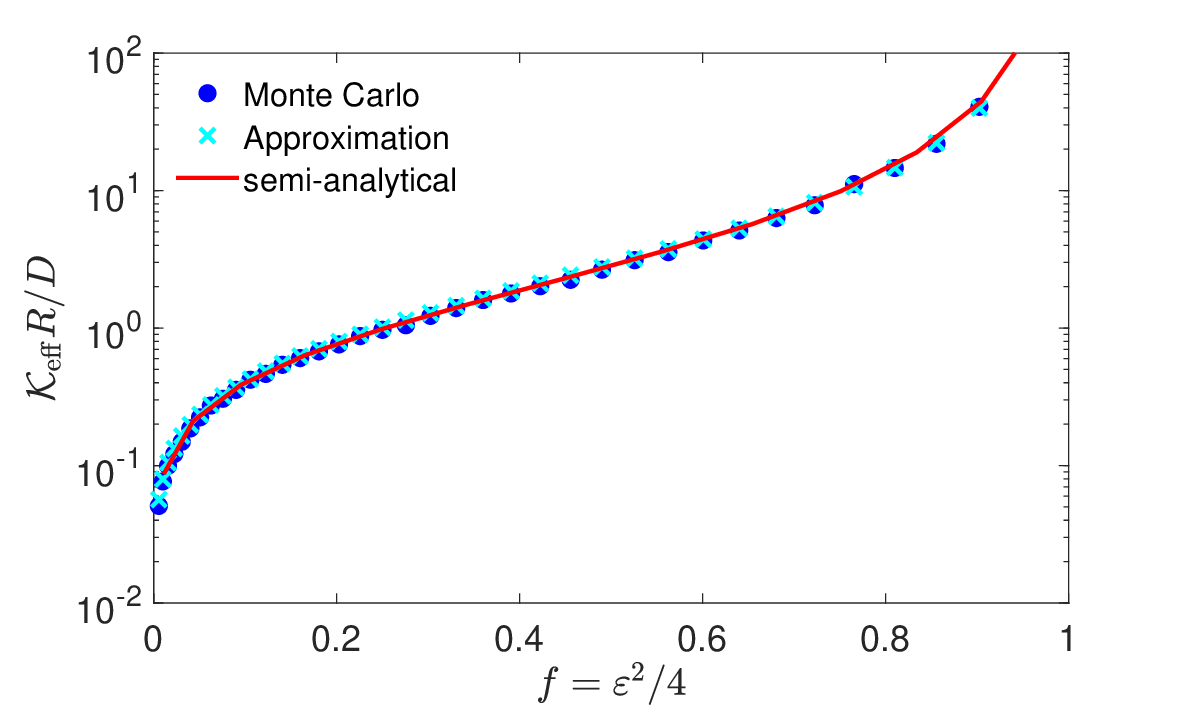} 
\end{center}
\caption{ Rescaled effective reactivity $\K_{\rm eff} R/D$ as a
  function of the patch surface fraction $f = \eps^2/4$ for a single
  absorbing patch of radius $\eps$ on the unit sphere (with
  $\K = \infty$).  Comparison between Monte Carlo simulations (with
  $M = 10^5$ realizations and $\ell = 10^{-2}$), the
  empirical approximation (\ref{eq:keff_Dagdug}), and semi-analytical
  solution from \cite{Traytak95}. }
\label{fig:Dagdug}
\end{figure}

(ii) For a single perfectly reactive circular patch of radius $\eps R$
and infinite reactivity $\K = \infty$, the mixed BVP (\ref{mfpt:ssp0})
can be reduced to dual infinite series relations that offers an
efficient semi-analytical solution.  In this way, Traytak
\cite{Traytak95} calculated the probability of reaction $P_{\rm T}$,
as given in the last column of Table 1 from \cite{Traytak95}, and
compared its values for different patch sizes $\eps$ with those
obtained by other numerical methods.  We used his semi-analytical
values to validate our Monte Carlo algorithm for perfectly reactive
patches.  Figure~\ref{fig:Dagdug} shows an excellent agreement between
our estimation of the effective reactivity $\K_{\rm eff}$ and that
derived from Traytak's solution, given by $\tfrac{R}{D} \K_{\rm T} =
(1/P_{\rm T} - 1)^{-1}$, over the whole range $0<f<1$ of patch area
fraction $f = \pi \eps^2 R^2/(4\pi R^2) = \eps^2/4$.  Moreover, we
also compare our results to an empirical approximation for the
effective reactivity given in \cite{Grebenkov25c} by
\begin{equation}  \label{eq:keff_Dagdug}  
  \frac{R}{D} \K_{\rm G} = \frac{2}{\pi} \sqrt{f} \, \frac{
      1 + 2.32\sqrt{f} - 1.47 f^2}{(1-f)^{3/2}} \,.
\end{equation}
(note that the original approximation from \cite{Dagdug16}
failed to describe the reactivity in the limit $f\to 1$).  We observe
from Fig.~\ref{fig:Dagdug} that our Monte Carlo numerical results
agree very closely with both Traytak's semi-analytical solution and
the empirical formula (\ref{eq:keff_Dagdug}) over the entire
range of $f$.

\begin{figure}
    \centering
     \begin{subfigure}[b]{0.49\textwidth}  
      \includegraphics[width =\textwidth]{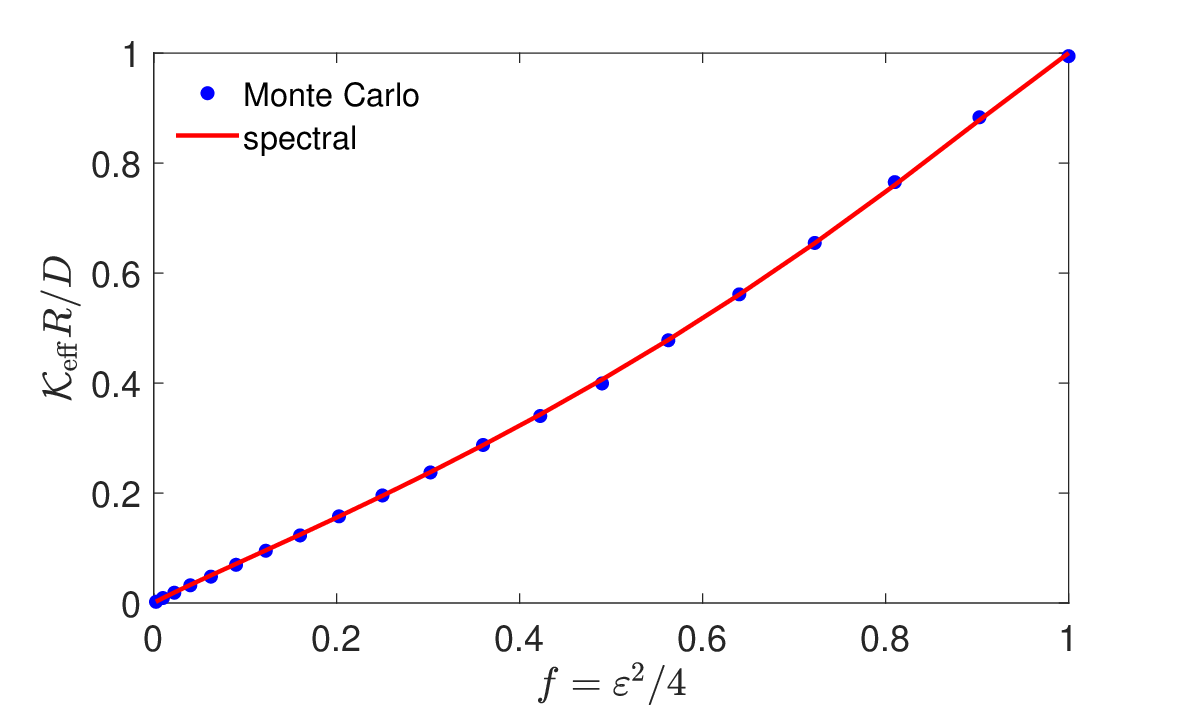} 
        \caption{$N = 1$}
        \label{fig:spectral_a}
    \end{subfigure}
    \begin{subfigure}[b]{0.49\textwidth}
      \includegraphics[width=\textwidth]{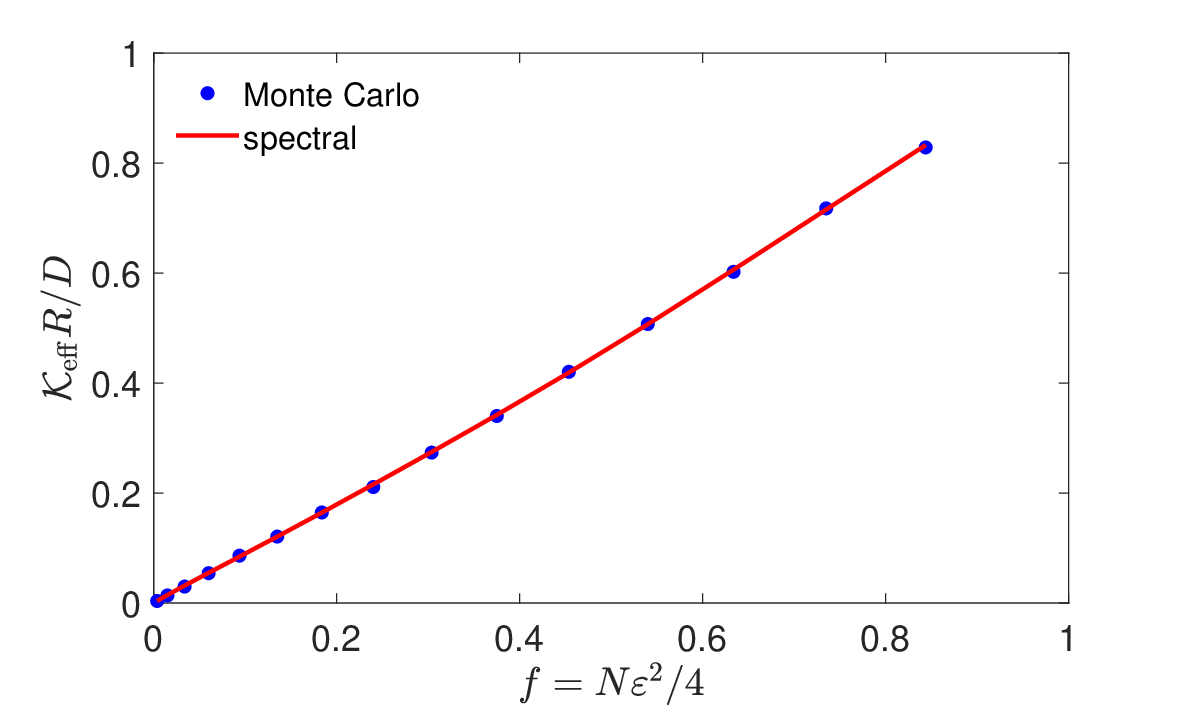} 
        \caption{$N = 6$} 
        \label{fig:spectral_b}
    \end{subfigure}
    \caption{Rescaled effective reactivity $\K_{\rm eff} R/D$ as a
      function of the patch surface fraction $f = N \eps^2/4$ for $N$
      identical partially reactive patches of radius $\eps$ on the
      unit sphere (with $\K R/D = 1$).  Comparison between Monte Carlo
      simulations (with $M = 10^5$ realizations and
      $\ell = 10^{-2}$) and the spectral solution from
      \cite{Grebenkov19b}, in which the matrices $\MM$ and $\KK$ were
      truncated to the size $121\times 121$.  (a) $N = 1$, $\eps$
      ranges from $0$ to $2$; (b) $N = 6$, with patches centered at the
      vertices of an octahedron, and $\eps$ ranges from $0$ to $0.75$.}
\label{fig:spectral}
%
\end{figure}

(iii) The case of partially reactive patches is least documented.  For
the validation of our algorithm, we compared our Monte Carlo results
with those computed using the spectral approach presented in
\cite{Grebenkov19b}.  This spectral method allows one to accurately
compute the flux $J$ onto a sphere with a finite number of circular
patches as $J = J_{\rm Smol} \, h_{00}^{(0)}$ (see Eq.~(52) of
\cite{Grebenkov19b}), where
$h_{00}^{(0)} = [(\MM + \KK)^{-1} \KK]_{00,00}$ is precisely the
probability of reaction. Here $\MM$ and $\KK$ are two
infinite-dimensional matrices that represent the Dirichlet-to-Neumann
operator and the reactivity distribution, respectively.  This
representation is exact and valid for any bounded distribution of
reactivity.  A practical implementation of this representation
requires truncating the matrices $\MM$ and $\KK$, whose efficient
construction for a sphere was provided in \cite{Grebenkov19b}.  Note
that higher truncation orders are needed to deal with larger
reactivities and/or smaller patches.  We used this numerical approach
to validate our Monte Carlo simulations for $N$ circular patches of
size $\eps$ with reactivity $R\K/D = 1$.  Fig.~\ref{fig:spectral_a}
presents the effective reactivity ${\mathcal K}_{\rm eff}$ of a single
patch as a function of its surface fraction $f$.  As expected, the
effective reactivity tends to $1$ as $f \to 1$, i.e., when the whole
surface is covered by the reactive patch.  In turn,
Fig. \ref{fig:spectral_b} illustrates the effective reactivity of six
identical partially reactive patches centered at the vertices of an
octahedron.  To avoid overlapping between patches, their common radius
$\eps$ is limited to $2\sin(\pi/8)\approx 0.7654$ that yields the
maximal covered fraction $0.88$.  In both cases, we observe an
excellent agreement between Monte Carlo simulations and the spectral
method, which provides a further benchmark for the accuracy of our
Monte Carlo algorithm.

\section{Numerical Comparison}\label{sec:comp}

We now employ the Monte Carlo algorithm of \S \ref{sec:numer} to
assess the accuracy and range of validity of our asymptotic results.
We recall that the homogenized formula (\ref{eq:keff}) (or
equivalently (\ref{homo:keff})) for the effective reactivity was
derived under the assumption of a large $N\gg 1$ number of small,
uniformly distributed, circular patches of radius $\eps$, in the 
limit of small patch area fraction: $f={\eps^2 N/4}\ll 1$.
Therefore, our numerical experiments will seek to explore the
following issues: (i) whether (\ref{eq:keff}) remains accurate for
moderately large $\eps$ and moderately small $N$; (ii) does
(\ref{eq:keff}) hold on the full range $\kappa>0$ of reactivity, and
(iii) how relevant is the assumption of equally spaced patches?

\begin{figure}
    \centering
     \begin{subfigure}[b]{0.49\textwidth}  
      \includegraphics[width =\textwidth]{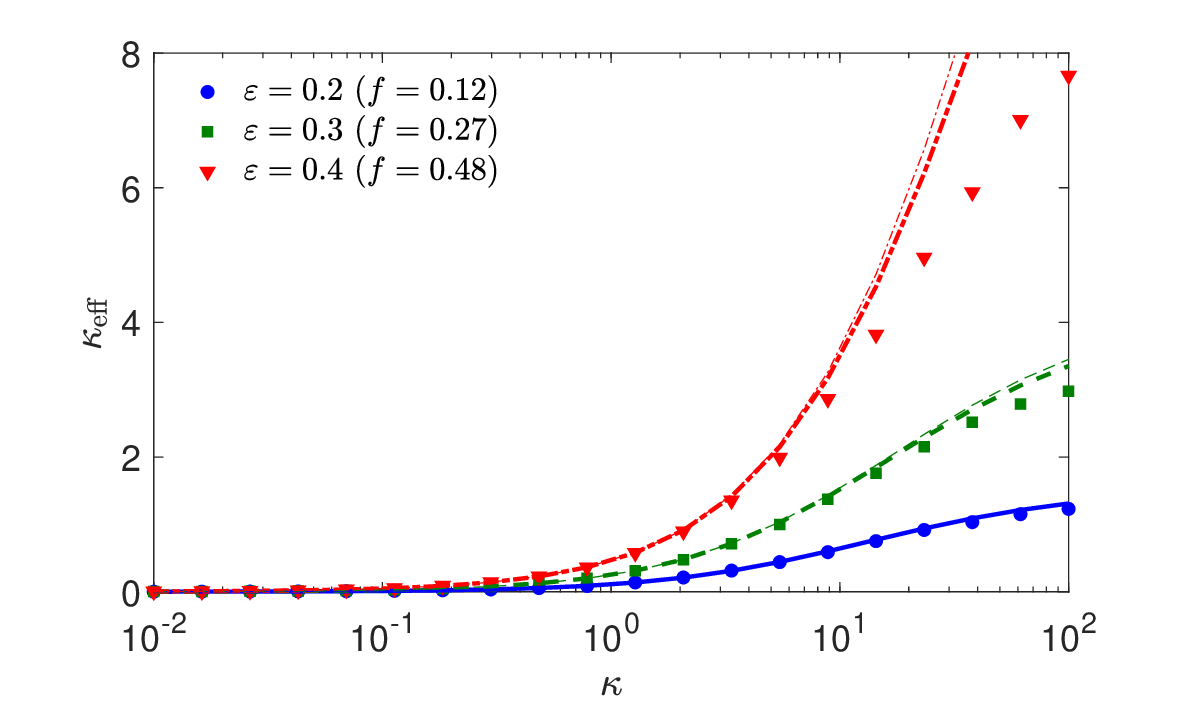} 
        \caption{$\keff$ versus $\kappa$}
        \label{fig:keff_N12a}
    \end{subfigure}
    \begin{subfigure}[b]{0.49\textwidth}
      \includegraphics[width=\textwidth]{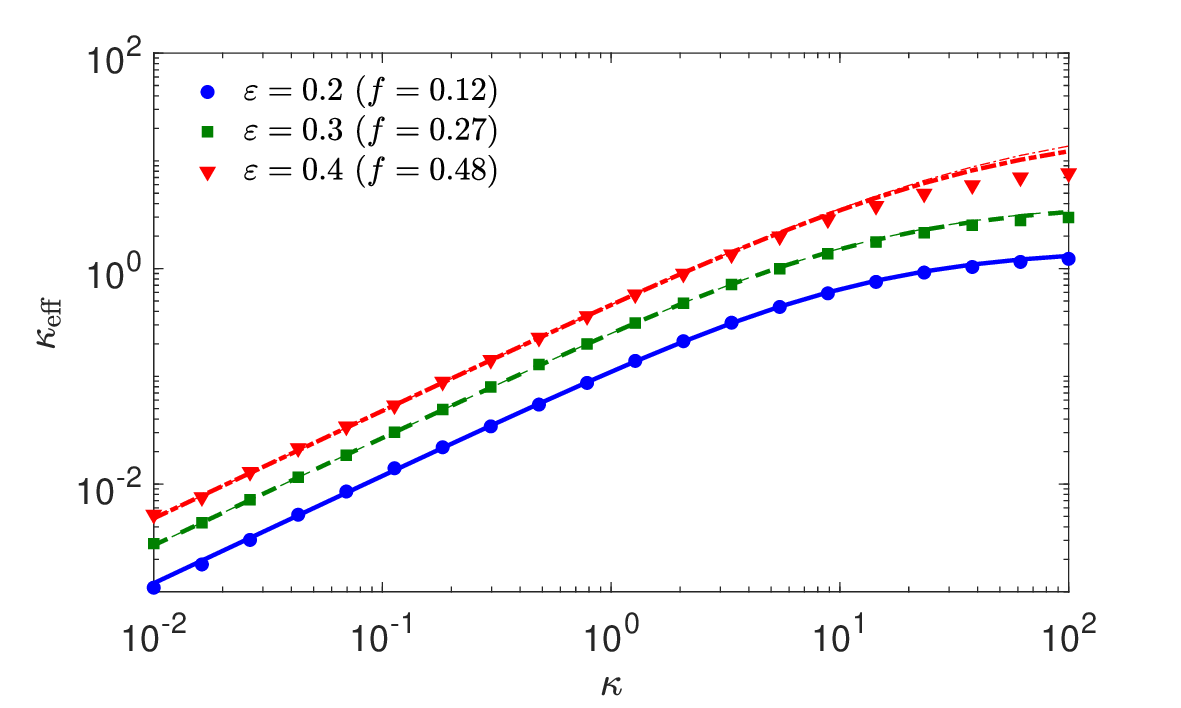} 
        \caption{$\keff$ versus $\kappa$ (log scale)} 
        \label{fig:keff_N12b}
    \end{subfigure}
\caption{(a): Dimensionless effective reactivity $\keff$ versus
  $\kappa = \eps \K R/D$ for $N = 12$ identical circular patches of
  radius $\eps$ centered at the vertices of an icosahedron.  Symbols
  are the Monte Carlo results (with $M = 10^5$ realizations and
  $\ell = 10^{-2}$), thick lines show (\ref{eq:kappa_J}) with $C_{\rm
  T}$ from the asymptotic formula (\ref{homo:c0_orig_1}), and thin
  lines are the homogenized asymptotic formula (\ref{eq:keff}).  (b):
  Same plot but with a logarithmic scale on the vertical axis to
  better show the small $\kappa$ comparison.}
\label{fig:keff_N12}
\end{figure}

\begin{figure}
    \centering
     \begin{subfigure}[b]{0.49\textwidth}  
      \includegraphics[width =\textwidth]{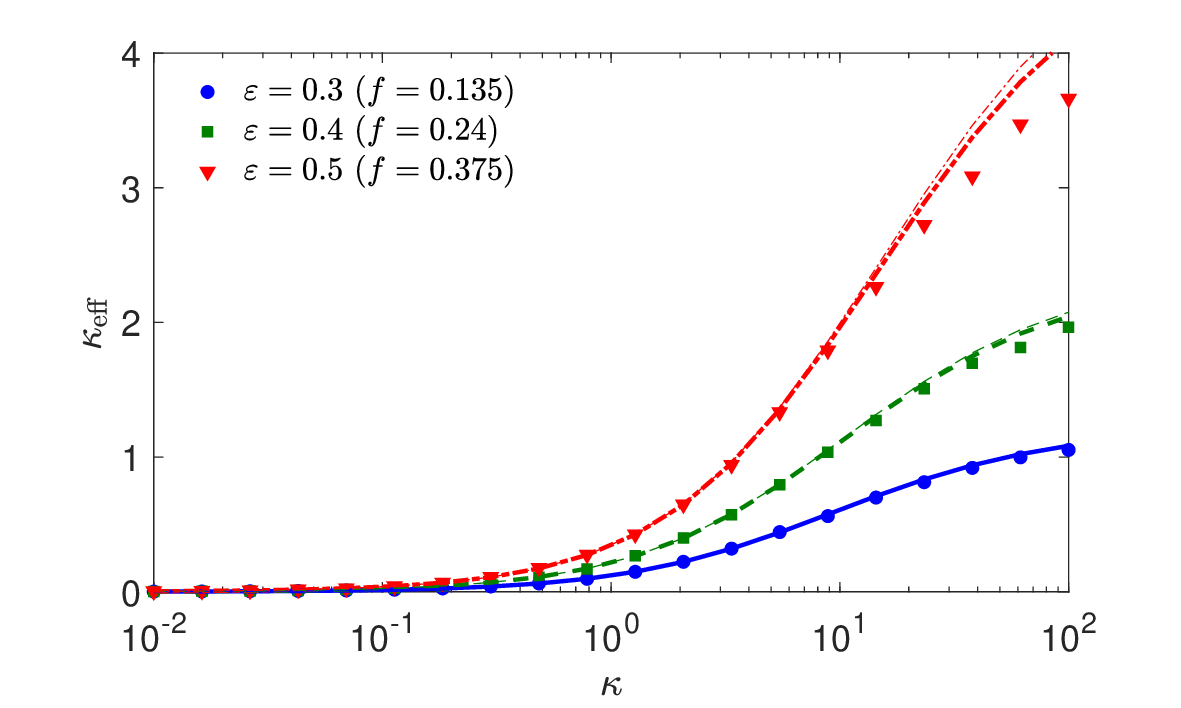} 
        \caption{$\keff$ versus $\kappa$}
        \label{fig:keff_N6a}
    \end{subfigure}
    \begin{subfigure}[b]{0.49\textwidth}
      \includegraphics[width=\textwidth]{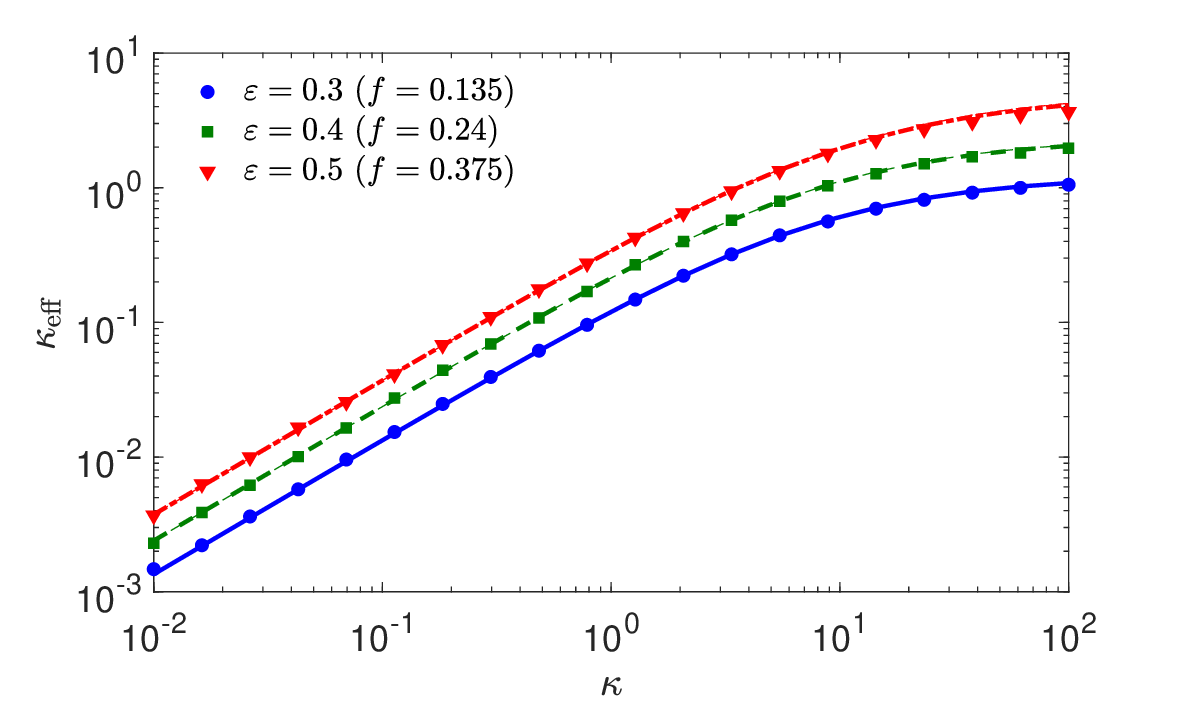} 
        \caption{$\keff$ versus $\kappa$ (log scale)} 
        \label{fig:keff_N6b}
    \end{subfigure}
\caption{(a): Dimensionless effective reactivity $\keff$ versus
  $\kappa = \eps \K R/D$ for $N = 6$ identical circular patches of
  radius $\eps$ centered at the vertices of an octahedron.  Symbols
  are the Monte Carlo results (with $M = 10^5$ realizations and
  $\ell =10^{-2}$), thick lines show (\ref{eq:kappa_J}) with $C_{\rm
  T}$ from the asymptotic formula (\ref{homo:c0_orig_1}), and thin
  lines are the homogenized asymptotic formula (\ref{eq:keff}). (b):
  Same plot but with a logarithmic scale on the vertical axis to
  better show the small $\kappa$ comparison.}
\label{fig:keff_N6}
\end{figure}

In Fig.~\ref{fig:keff_N12} we plot the effective reactivity
$\kappa_{\rm eff}$ as a function of $\kappa$ for 12 patches centered
at the vertices of an icosahedron on the unit sphere, and with three
choices of the patch radius $\eps$, that give different patch coverage
fractions $f$.  We compare Monte Carlo results (shown by symbols) to
two asymptotic formulas: the general formula (\ref{eq:kappa_J}), in
which the capacitance $C_{\rm T}$ is given in (\ref{homo:c0_orig_1}), and
the homogenized formula (\ref{eq:keff}).  For patches of moderate size
$\eps = 0.2$ that cover one-eighth of the spherical surface, both
asymptotic formulas are in remarkable agreement with Monte Carlo
simulations over the whole considered range of reactivities from
$10^{-2}$ to $10^2$.  When $\eps = 0.3$, corresponding to a
one-quarter coverage, the agreement is still excellent for small and
moderate reactivities, but small deviations appear at high
reactivities.  These deviations at large $\kappa$ are further enhanced
when $\eps = 0.4$, which corresponds to a one-half coverage of the
sphere by patches.  Despite these deviations, the asymptotic formula
is still applicable for small and moderate reactivities.  The origin
of this remarkable agreement is revealed by the structure of the
asymptotic formula (\ref{eq:keff}), in which the patch size $\eps$ is
multiplied by the reactive capacitance $C(\kappa)$ (except for
the logarithmic term).  In other words, the actual small parameter is
$\eps C(\kappa)$, which vanishes in the limit $\kappa \to 0$ since
$C(\kappa)\sim{\kappa/2}$ for $\kappa\ll 1$. As shown in
Fig.~\ref{fig:keff_N6}, similar results are obtained for $N = 6$
patches that are centered at the vertices of an octahedron.

\begin{figure}
\begin{center}
\includegraphics[width=88mm]{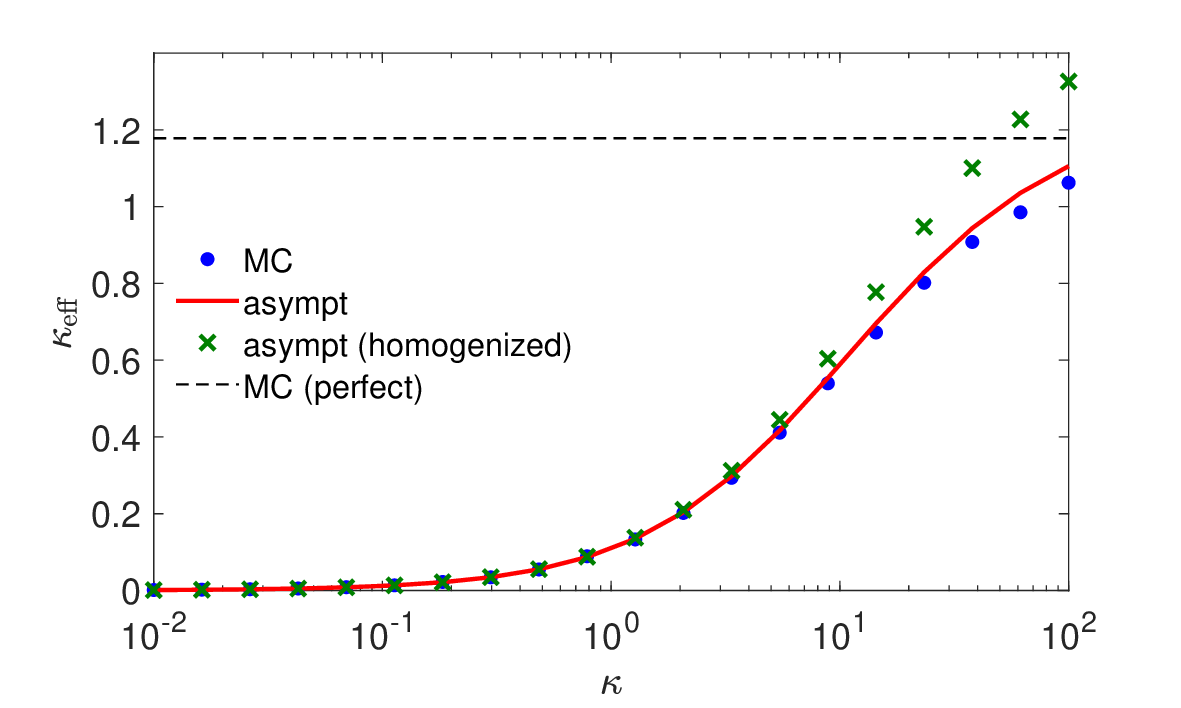} 
\end{center}
\caption{ Dimensionless effective reactivity $\keff$ versus
  $\kappa = \eps \K R/D$ for $N = 12$ identical circular patches of
  radius $\eps = 0.2$, that are centered at independent random points
  uniformly distributed on the sphere.  Round symbols are the Monte
  Carlo results (with $M = 10^5$ realizations and
  $\ell = 10^{-2}$), the thick line shows (\ref{eq:kappa_J})
  with $C_{\rm T}$ from the asymptotic formula (\ref{homo:c0_orig_1}),
  and crosses are the homogenized asymptotic formula (\ref{eq:keff}).
  The thin horizontal line is the Monte Carlo computed value for
  a perfectly reactive patch.}
\label{fig:keff_N12_rand}
\end{figure}

In Fig.~\ref{fig:keff_N12_rand} we illustrate the effect of random
locations of twelve identical circular patches of radius $\eps = 0.2$.
The general asymptotic formula (\ref{eq:kappa_J}) with
(\ref{homo:c0_orig_1}) yields accurate results over the whole
considered range of reactivities.  As the patches are not equally
spaced, the homogenized asymptotic formula (\ref{eq:keff}) is
expectedly less accurate at high reactivities, but still provides a
reasonable approximation at smaller reactivities.

\begin{figure}
\begin{center}
\includegraphics[width=88mm]{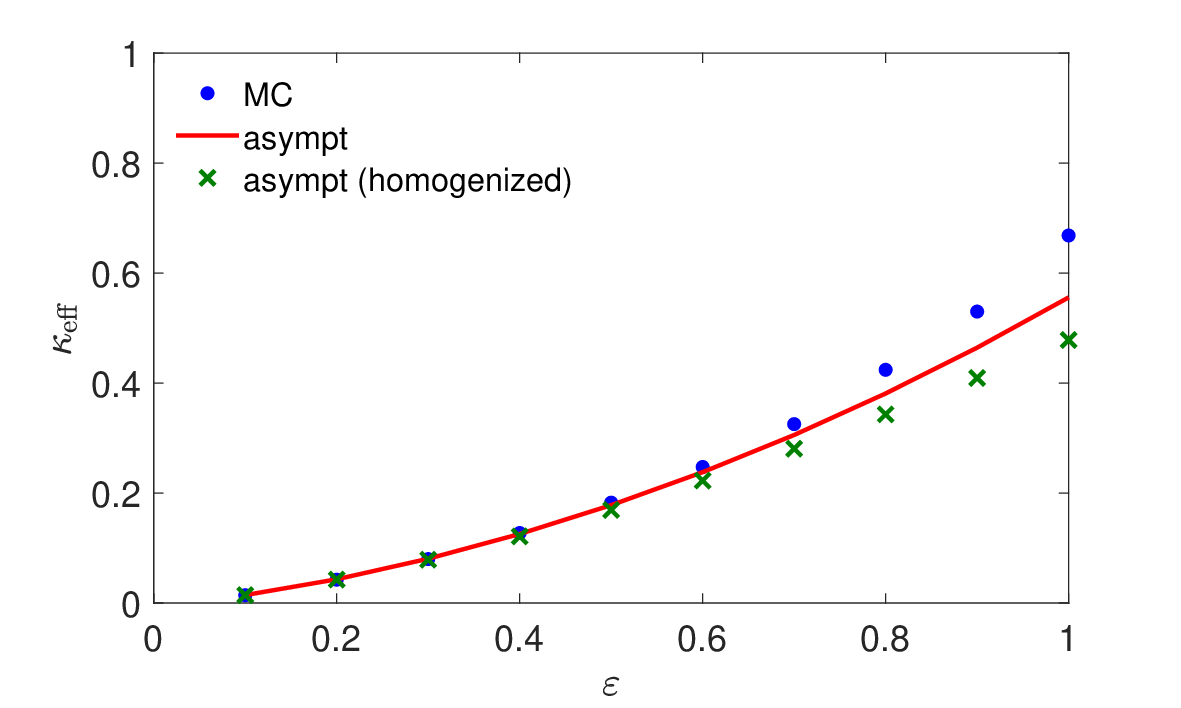} 
\end{center}
\caption{Dimensionless effective reactivity $\keff$ versus $\eps$ for
  one circular patch of radius $\eps$ with reactivity $R\K/D = 10$.
  Round symbols are Monte Carlo results (with $M = 10^5$
  realizations and $\ell = 10^{-2}$), the thick line shows
  (\ref{eq:kappa_J}) with $C_{\rm T}$ from the asymptotic formula
  (\ref{homo:c0_orig_1}), and crosses are the homogenized asymptotic
  formula (\ref{eq:keff}). The patch area fraction $f={\eps^2/4}$ is
  $25\%$ when $\eps=1$.}
\label{fig:keff_N1}
\end{figure}

Finally, in Fig.~\ref{fig:keff_N1} we consider the extreme case $N =
1$ of a single reactive patch and we plot the effective reactivity as
a function of the patch radius $\eps$, for a fixed moderately large
reactivity $R\K/D = 10$.  Even for large patches for which $\eps = 1$,
corresponding to a one-quarter coverage of the sphere, the general
asymptotic formula (\ref{eq:kappa_J}) with (\ref{homo:c0_orig_1}), as
well as the homogenized asymptotic formula (\ref{eq:keff}), remain
reasonably accurate.  Their accuracy is further improved when
$\K$ is decreased.

\section{Discussion}\label{sec:discussion}

By using an asymptotic approach based on strong localized perturbation
theory \cite{Ward93} we have derived a three-term asymptotic formula
for the capacitance of a spherical target that has many small surface
patches of finite reactivity on an otherwise reflecting boundary. This
problem has broad applications in biophysical modeling, where a common
theme in diverse applications is the study of how diffusing ligands
can bind to surface receptors that have finite reactivity. Our
analysis, relying on a non-conventional geodesic normal
coordinate system, has extended the previous analysis in
\cite{Lindsay17} for perfectly reactive circular patches on a sphere
and the leading-order analysis in \cite{Plunkett24} for partially
reactive circular patches on an infinite plane. Our
  three-term asymptotic result for the capacitance in the small patch
  radius limit is identical in form to that derived in
  \cite{Lindsay17}, provided that we now use the reactive capacitance
  of a patch and re-define an additional monopole coefficient. This
  three-term result incorporates the effect of inter-patch
  interactions, accounts for the curvature of the boundary of the
  sphere, and remains well-ordered over the full range of
  reactivities. By homogenizing our result for the capacitance we
have derived a scaling law for the effective reactivity of the
structured target that has identical, uniformly distributed, partially
reactive small patches. The two coefficients, depending on the local
reactivity and the shape of the patches, that appear in our scaling
law can be calculated from Steklov eigenfunction expansions, while for
circular patches they are well-approximated by heuristic
approximations.  Finally, our asymptotic results have been validated
for certain patch configurations by numerical results obtained from a
new Monte Carlo algorithm. This comparison has suggested that the
asymptotic theory can still be reliably used beyond its expected range
of validity (small patches).

We remark that it is readily possible to extend the leading-order
analysis in \cite{Plunkett24} by deriving a three-term expansion for
the capacitance of an infinite plane that has small circular partially
reactive patches on an otherwise reflecting boundary. For this simpler
scenario, where boundary curvature plays no role, there is no monopole
coefficient $E$ and no subdominant logarithmic singularity of the
surface Green's function. A similar analysis for circular patches of
infinite reactivity, with either random or periodic patch locations on
the infinite plane, was undertaken in \cite{Lindsay18a} and
\cite{Lindsay18b}, respectively. More specifically, consider a
collection of $N$ well-separated disks centered at $\x_i$, for
$i=1,\ldots,N$, on the infinite plane, each with radius $\eps a_i$ and
reactivity $\kappa_i$. Then, by a simple adaptation of the analysis in
\S 2 of \cite{Lindsay18a}, which dealt with the case of
infinite reactivity, we obtain that the asymptotic expansion for the
dimensionless capacitance $C_{\rm T}$ of this planar structured target
is
\begin{equation}\label{disc:cap}
    \frac{1}{C_{\rm T}}\sim \frac{1}{\eps \overline{C}} \left[
      1 + \eps \gamma + \eps^2\left(\gamma^2 -\frac{1}{\overline{C}}
        \sumcapone\sumcaptwo\sumcapthree
        \frac{C_i C_j C_k}{|\x_k-\x_j||\x_j-\x_i|}
        \right)\right]\,,
\end{equation}
where $\gamma\equiv {\vc^{T}{\mathcal G}_s\vc/\overline{C}}$,
$\vc\equiv (C_1,\ldots,C_N)^T$,
$\overline{C}\equiv\overline\sum_{j=1}^{N}C_j$, and
$C_j=a_j{\mathcal C}(a_j \kappa_j)$ is the reactive capacitance
defined by the inner problem (\ref{mfpt:wc}). Here the $N\times N$
symmetric Green's matrix ${\mathcal G}_s$ needed for $\gamma$ has
matrix entries $\left({\mathcal G}_s\right)_{ii}=0$ and
$\left({\mathcal G}_s\right)_{ij}=|\x_i-\x_j|^{-1}$ for $i\neq j$. We
remark that this result can be readily obtained by simply replacing
the capacitance of a perfectly absorbing disk appearing in
\cite{Lindsay18a} with the reactive capacitance of a patch. In this
way, by using the heuristic approximation (\ref{mfpt:sigmoidal_2}),
(\ref{disc:cap}) is readily evaluated over the full range of
reactivities of the patches. Moreover, the analysis in
\cite{Lindsay18b} that determined an effective capacitance for
periodic patterns of identical circular patches centered on a Bravais
lattice on the infinite plane with unit area of the fundamental
Wigner-Seitz cell is easily extended to the case of partial
reactivity.

Finally, we highlight two problems that warrant further study. Firstly
by exploiting a local tangential-normal coordinate system together
with a careful resolution of the subdominant logarithmic singularity
of the surface Green's function (see \cite{Efimov74,Silber03}), it
should be possible to derive a similar three-term result for the
capacitance of a generic bounded domain with a smooth boundary covered
by small surface patches of finite reactivity. Such a result would
also depend on the mean curvature of the boundary at the center of
each patch as well as the regular (i.e., non-singular) part of the
surface Green's function at each patch location, the latter of which
must be computed numerically. With these modifications, the overall
analysis should be rather similar to that done for the structured
sphere in \S \ref{mfpt_sec:expan} and the Appendices. Secondly, from a
computational viewpoint, it would be worthwhile to devise a numerical
PDE approach to accurately solve (\ref{bp:ssp}) for a large collection
of partially reactive patches. A key challenge in the numerics is to
carefully resolve the behavior at the boundary of each patch. For
perfectly reactive patches, corresponding to the Dirichlet-Neumann
problem, such a scheme has previously been developed
(cf. \cite{Lindsay17,Lindsay18a,Kaye20}).

\section*{Acknowledgements} 
D.S.G. acknowledges the Simons Foundation for supporting his
sabbatical sojourn in 2024 at the CRM (CNRS -- University of
Montr\'eal, Canada), and the Alexander von Humboldt Foundation for
support within a Bessel Prize award.  M.J.W. was supported by the
NSERC Discovery grant program. We are grateful to Prof. Sean Lawley of
the University of Utah for discussions related to \cite{Plunkett24}
and to the sigmoidal approximation in (\ref{mfpt:sigmoidal_2}).

\appendix
\renewcommand{\theequation}{\Alph{section}.\arabic{equation}}

\section{Geodesic Normal Coordinates to the Unit Sphere
$\Omega$}\label{app_g:geod}

We define geodesic normal coordinates
$\bxi=(\xi_1,\xi_2,\xi_3)^T\in \left({-\pi/2},{\pi/2}\right) \times
\left(-\pi,\pi\right)\times [0,\infty]$ in the exterior $|\x|\geq 1$
of the unit sphere $\Omega$. In these coordinates, $\bxi=0$
corresponds to $\x_i\in\partial\Omega$, while $\xi_3>0$ holds in the
exterior of $\Omega$.  In terms of the spherical angles
$\theta_i\in (0,\pi)$ and $\varphi_i\in [0,2\pi)$, and for $|\x_i|=1$, we
first define the orthonormal vectors $\x_i$, $\vv_{2i}$ and $\vv_{3i}$ by
\begin{equation}\label{app_g:ortho}
  \x_i \equiv \begin{bmatrix}
           \cos\varphi_i \, \sin\theta_i \\
           \sin\varphi_i \, \sin\theta_i \\
           \cos\theta_i
         \end{bmatrix}\,, \,\,\,
         \vv_{2i}=\partial_{\theta} \x_i \equiv
\begin{bmatrix}
           \cos\varphi_i \, \cos\theta_i \\
           \sin\varphi_i \, \cos\theta_i \\
           -\sin\theta_i
         \end{bmatrix}\,, \,\,\,
         \vv_{3i}=\x_i {\bf \times} \partial_{\theta} \x_i \equiv
\begin{bmatrix}
           -\sin\varphi_i \\
           \cos\varphi_i \\
           0
         \end{bmatrix}\,,
\end{equation}
where $\vv_{2i}$ and $\vv_{3i}$ provide a basis for the tangent plane
to the sphere at $\x=\x_i$. In terms of these vectors, the geodesic
coordinates $\bxi=(\xi_1,\xi_2,\xi_3)^T$ are defined by the global
transformation
\begin{equation}\label{app_g:global}
  \x(\bxi) = \left(1+\xi_3\right)\left( \cos\xi_1 \,
    \cos\xi_2 \, \x_i + \cos\xi_1 \,
    \sin\xi_2 \, \vv_{2i} + \sin\xi_1 \vv_{3i}\right)\,.
\end{equation}
The geodesic coordinate curves obtained by setting $\xi_3=0$, and
fixing either $\xi_2=0$ or $\xi_1=0$ are, respectively,
$\x(\xi_1,0,0)=\cos\xi_1 \, \x_i + \sin\xi_1 \, \vv_{3i}$ or
$\x(0,\xi_2,0)=\cos\xi_2 \, \x_i + \sin\xi_2 \, \vv_{2i}$. These
circles correspond to intersections of
$\partial\Omega$ with planes spanned by
$\lbrace{\x_i,\vv_{3i}\rbrace}$ or $\lbrace{\x_i,\vv_{2i}\rbrace}$,
respectively.

The scale factors $h_{\xi_j}\equiv \vert {\partial\x/\partial \xi_j}\vert$
for $j=1,2,3$ are readily calculated as
\begin{equation}\label{app_g:scale}
  h_{\xi_1} = (1+\xi_3) \,, \qquad h_{\xi_2}=(1+\xi_3)\cos\xi_1 \,, \qquad
  h_{\xi_3}=1 \,.
\end{equation}
For a generic function ${\mathcal V}(\bxi)\equiv
u\left(\x({\bxi})\right)$, we calculate, as similar to that done in
Appendix A of \cite{GrebWard25}, that the Laplacian transforms
according to
\begin{equation}\label{app_g:lap}
\Delta_{\x} u = {\mathcal V}_{\xi_3\xi_3} + \frac{2}{1+\xi_3} {\mathcal V}_{\xi_3} +
\frac{1}{(1+\xi_3)^2\cos^{2}\xi_1} {\mathcal V}_{\xi_2\xi_2} +
\frac{1}{(1+\xi_3)^2 \cos\xi_1} \frac{\partial}{\partial_{\xi_1}}
\left(\cos\xi_1 {\mathcal V}_{\xi_1} \right)\,.
\end{equation}

Then, upon introducing the local (or inner) variables,
$\y=(y_1,y_2,y_3)^T$, defined by
\begin{equation}\label{app_g:innvar}
  \xi_1=\eps y_1 \,, \qquad \xi_2=\eps y_2 \,, \qquad \xi_3 =\eps y_3\,,
\end{equation}
we use the Taylor series approximations $(1+\xi_3)^{-1}\sim 1-\eps y_3$,
$(1+\xi_3)^{-2}\sim 1-2\eps y_3$, $\cos^{2}\xi_1=1+{\mathcal O}(\eps^2)$ and
$\sin\xi_1\sim \eps y_1$, to show that (\ref{app_g:lap}) reduces
to (\ref{mfpt:local}).

To determine a two-term approximation for the Euclidian distance
$|\x-\x_i|$ near the patch, we proceed in a similar way as in
Appendix A of \cite{GrebWard25}. By substituting (\ref{app_g:innvar}) in
(\ref{app_g:global}), we obtain from a Taylor series approximation that
\bsub \label{app_g:vab_all}
\begin{equation}
  \x-\x_i = \eps \vb_0 - \eps^2 \vb_1 + {\mathcal O}(\eps^3)\,, \qquad
  |\x-\x_i|^2 \sim \eps^2 \left( \vb_0^T\vb_0 -2\eps \vb_0^T\vb_1\right)\,,
\end{equation}
where $\vb_0$ and $\vb_1$ are defined by
\begin{equation}\label{app_g:vab}
  \vb_0 =  y_3 \x_i + y_2 \vv_{2i} + y_1 \vv_{3i} \,, \qquad
  \vb_1= \frac{1}{2} \left(y_1^2+y_2^2\right) \x_i - y_3 y_2 \vv_{2i} - y_3
  y_1 \vv_{3i} \,.
\end{equation}
\esub
In this way, and by labeling $\rho=|\y|$, with $\y=(y_1,y_2,y_3)^T$, we
conclude that
\begin{equation}\label{app_g:loc}
  |\x-\x_i|\sim \eps \rho + \frac{\eps^2 y_3}{2\rho} \left(y_1^2+y_2^2
    \right) + {\mathcal O}(\eps^3) \,, \quad
  \frac{1}{|\x-\x_i|} \sim \frac{1}{\eps \rho} \left( 1 -
    \frac{\eps y_3}{2\rho^2} (y_1^2+y_2^2) + {\mathcal O}(\eps^2)\right)\,.
\end{equation}

In matrix form, and to the leading order in $\eps$, we can write
(\ref{app_g:vab_all}) in terms of $\y=(y_1,y_2,y_3)^T$ and an
orthogonal matrix ${\mathcal Q}_i$ as
\begin{equation}\label{app_g:change}
  \y \sim \eps^{-1} {\mathcal Q}_{i}^T (\x-\x_i) \,, \quad \mbox{where} \quad
  {\mathcal Q}_{i} \equiv \begin{bmatrix}
    \vert & \vert & \vert \\
    \vv_{3i}   & \vv_{2i} & \x_i \\
    \vert & \vert & \vert
\end{bmatrix} \quad \rightarrow \quad |\y|\sim \eps^{-1}|\x-\x_i|\,.
\end{equation}
Finally, since $|\x-\x_i|=\eps \rho + {\mathcal O}(\eps^3)$ when
$y_3=0$ from (\ref{app_g:loc}), and by recalling the scale factor
$h_{\xi_3}=1$ from (\ref{app_g:scale}), we conclude that
a Robin boundary condition on a locally circular patch is well-approximated
in the local geodesic coordinates by
\begin{equation}
  - \partial_{y_3} U  + \kappa U=0 \,, \quad
  \mbox{for} \quad y_3=0 \,, \,\, (y_1^2+y_2^2)^{1/2}\leq a +
  {\mathcal O}(\eps^2)\,. \label{app_g:robin}
\end{equation}
To the order of our asymptotic expansion we can neglect the
${\mathcal O}(\eps^2)$ term in (\ref{app_g:robin}).

\section{Computing the Reactive Capacitance}
\label{app:Cmu}

In this Appendix, we summarize some of the results from Appendix D of
\cite{GrebWard25} that determined a Steklov eigenfunction expansion for
the reactive capacitance $C_i(\kappa_i)$ for an arbitrary patch shape
$\Gamma_i$, and for the special case of a circular patch.

Written in terms of the local geodesic coordinates, the following
Steklov eigenvalue problem for eigenpairs $\Psi_{ki},\mu_{ki}$ in a
half-space $\R_{+}^{3}$ plays a central role in determining
$C_i(\kappa_i)$:
\begin{subequations}  \label{eq:Psi_def}
\begin{align}  \label{eq:Vk_eq}
\Delta \Psi_{ki} & = 0 \,, \quad {\bf y} \in \R_+^3 \,,\\
  \partial_n \Psi_{ki} & = \mu_{ki} \Psi_{ki}\,, \quad y_3=0 \,,\,
                         (y_1,y_2)\in \PT_i\,,
  \\  \label{eq:Vk_Neumann}
  \partial_n \Psi_{ki} & = 0 \,, \quad y_3=0 \,,\, (y_1,y_2)\notin \PT_i\,,
  \\  \label{eq:Vk_inf}
  \Psi_{ki}(\y) & = {\mathcal O}\left({1/|\y|}\right) \,,
                  \quad \textrm{as}\quad |\y|\to \infty\,.
\end{align}
\end{subequations}
As discussed in Appendix D of \cite{GrebWard25}, the Steklov
eigenvalues, enumerated by the index $k = 0,1,2,\ldots$, can be
ordered as
\begin{equation}
  0 < \mu_{0i} < \mu_{1i} \leq \cdots \nearrow +\infty \,,
\end{equation}
where the principal eigenvalue $\mu_{0i}$ is simple and strictly
positive.  The corresponding eigenfunctions, when restricted to the
patch $\Gamma_i$, as labeled by $\Psi_{ki}(\y)|_{\Gamma_i}$, form a complete
orthonormal basis in $L^2(\Gamma_i)$, in the sense that
$\int\limits_{\Gamma_i} \Psi_{ki} \Psi_{k^{\prime}i} \, d\y=
\delta_{k,k^{\prime}}$.  By applying the divergence theorem to
(\ref{eq:Psi_def}), the far-field behavior of
$\Psi_{ki}(\y)$ has the form
\begin{equation}  \label{eq:Psi_asympt}
  \Psi_{ki}(\y) \sim \frac{\mu_{ki} d_{ki}}{2\pi |\y|} + \ldots \,, \quad
  \mbox{as}   \quad |\y|\to\infty\,, \quad \mbox{where} \quad
  d_{ki} = \int_{\Gamma_i} \Psi_{ki} \, d\y\,.
\end{equation}

In terms of these Steklov eigenpairs and weights $d_{ki}$, it was
shown in Appendix D of \cite{GrebWard25} that
\begin{equation} \label{eq:Cmu}
  C_{i}(\kappa_i) = \frac{\kappa_i}{2\pi} \sum\limits_{k=0}^\infty
  \frac{\mu_{ki}  d_{ki}^2}
  {\mu_{ki} + \kappa_i} \,. 
\end{equation}
From this Steklov eigenfunction expansion, we conclude that
$C_{i}(\kappa_i)$ is monotonically increasing on $\kappa_i>0$, so that
$C_{i}(\infty)$ is an upper bound for $C_{i}(\kappa_i)$ on $\kappa_i>0$.
Moreover, for $\kappa_i\to 0^{+}$, a Taylor series approximation yields that
\begin{equation} \label{eq:Cmu_Taylor}
  C_{i}(\kappa_i) = - a_i \sum\limits_{n=1}^{\infty} c_{ni} \,
  \left(-\kappa_i a_i\right)^n \,, \quad
  \textrm{with}\quad c_{ni} = \frac{1}{2\pi a_i^{n+1}} \sum\limits_{k=0}^\infty
  \frac{d_{ki}^2}{\mu_{ki}^{n-1}} \,,
\end{equation}
where $c_{1i}={|\Gamma_i|/(2\pi)}$ (see Appendix D of \cite{GrebWard25}).

\subsection{Circular Patch}\label{app:Cmu_disk}

For a circular patch $\Gamma_i$ of unit radius, the first eight
axially symmetric eigenpairs for (\ref{eq:Psi_def}) were computed
numerically in \cite{Grebenkov24} (see also Appendix D.2 of
\cite{GrebWard25}). The eigenvalues and corresponding weights are
given in Table \ref{table:muk_disk}.

To characterize the limiting asymptotics of $C_{i}(\kappa_i)$ for
$\kappa_i\ll 1$, in Appendix B of \cite{GrebWard25} it was shown
analytically that the first three Taylor coefficients in
(\ref{eq:Cmu_Taylor}) are
\begin{equation}  \label{eq:cn_exact}
  c_{1i} = \frac12 \,,  \quad c_{2i} = \frac{4}{3\pi} \approx 0.4244\,,  \quad
  c_{3i} = \frac{4}{\pi^2} \int_{0}^{1} r \left[E(r)\right]^2
\, dr \approx 0.3651\,, 
\end{equation}
where $E(r)$ is the complete elliptic integral of the second
kind. Moreover, Appendix D.2 of \cite{GrebWard25} established that
all of the Taylor coefficients $c_{ni}$ in (\ref{eq:Cmu_Taylor}) are
well-approximated by
\begin{equation}  \label{eq:cn_approx}
  c_{ni} \approx \frac{0.4888}{(1.1578)^{n-1}} +
  \frac{0.0084}{(4.3168)^{n-1}}\,,
  \quad \mbox{for} \quad n\geq 2\,.
\end{equation}

\begin{table}
\begin{center}
\begin{tabular}{|c|c|c|c|c|c|c|c|c|} \hline
$k$     & 0      & 1      & 2      & 3      & 4      & 5      & 6      & 7      \\ \hline
$\mu_{ki}$ & 1.1578 & 4.3168 & 7.4602 & 10.602 & 13.744 & 16.886 & 20.028 & 23.169 \\ 
  $d_{ki}$   & 1.7524 & 0.2298 & 0.1000 & 0.0587 & 0.0397 & 0.0291 & 0.0225 & 0.0180\\
\hline
\end{tabular}
\end{center}
\caption{ The first eight Steklov eigenvalues $\mu_{ki}$ and weights
  $d_{ki}$ for the unit disk $\Gamma_i$ in the upper half-space that correspond
  to axially symmetric eigenfunctions on the patch, for which $d_{ki}\neq 0$,
  as computed numerically
  in \cite{Grebenkov24}.}
\label{table:muk_disk}
\end{table}

In contrast, in the limit $\kappa_i\to +\infty$, it was shown in
\cite{Guerin23} that the difference $C_{i}(\kappa_i)-C_i(\infty)$ is
not analytic in $\kappa_i$ for $\kappa_i\gg 1$. In particular, the
results in \cite{Guerin23} (see also Appendix D.3 of
\cite{GrebWard25}) yield the refined asymptotic behavior given in
(\ref{eq:Cmu_asympt}).

\section{Inner Problem Beyond Tangent Plane Approximation}\label{app_h:inn2}

In this Appendix we construct the solution to (\ref{mfpt_b:Phi2}) and
show how to determine the monopole coefficient $E_i$ in the limiting
behavior (\ref{mfpt_b:Phi2_ff}). Since this Appendix is similar to
Appendix C of \cite{GrebWard25} for the narrow capture MFPT problem
{\em inside} the unit sphere, we only briefly outline the
analysis. The key distinction from the analysis in \cite{GrebWard25}
is that for our {\em exterior} problem one must account for the
different algebraic sign of the curvature of the sphere, as viewed
from the exterior of the sphere.

The central issue in solving (\ref{mfpt_b:Phi2}) for $\Phi_{2i}$ is to
find an explicit particular solution $\Phi_{2pi}$ that accounts for
the inhomogeneous term in the PDE (\ref{mfpt_b:Phi2_1}) for
$\Phi_{2i}$. This inhomogeneous term is directly responsible for the
non-monopole behavior in the far-field (\ref{mfpt_b:Phi2_4}). However,
a second issue is that we must also account for the fact that this
particular solution $\Phi_{2pi}$ does not satisfy the Robin boundary
condition (\ref{mfpt_b:Phi2_2}) on the patch. As a result, in our
decomposition of $\Phi_{2i}$ we need to introduce an auxiliary
function $\Phi_{2hi}$, which satisfies the homogeneous part of the PDE
(\ref{mfpt_b:Phi2_1}), but that allows the homogeneous Robin condition
(\ref{mfpt_b:Phi2_2}) for the full solution $\Phi_{2i}$ to hold. The
far-field behavior of this auxiliary function yields the monopole
coefficient $E_i$. Our decomposition is summarized as follows:

\begin{lemma}\label{lemma:Phi2}
  The solution to (\ref{mfpt_b:Phi2}) can be decomposed as
\begin{equation}\label{app_h:decomp}
  \Phi_{2i} = \Phi_{2pi} + \Phi_{2hi} \,,
\end{equation}
where 
\begin{equation}\label{app_h:phi_2p}
  \Phi_{2pi}=\frac{y_3^2}{2}w_{i,y_3} +\frac{y_3}{2}w_i - \frac{1}{2}
  \int_{0}^{y_3} w_{i}(y_1,y_2,\eta;\kappa_i)\, d\eta -
  {\mathcal F}_{i}(y_1,y_2;\kappa_i)\,,
\end{equation}
with $w_i$ being the solution to (\ref{mfpt:wc}). Here
${\mathcal F}_{i}(y_1,y_2;\kappa_i)$, with $\Delta_{S}{\mathcal F}_{i} \equiv
{\mathcal F}_{i,y_1 y_1}+ {\mathcal F}_{i,y_2 y_2}$, is the unique solution to
(\ref{mfpt:fprob}), while $\Phi_{2hi}$ is the unique solution to
(\ref{mfpt:inn2_probh}). For an arbitrary patch shape $\PT_i$, the monopole
coefficient $E_i=E_i(\kappa_i)$ in (\ref{mfpt:inn2_h4}) is given by
\bsub \label{app_h:eval}
\begin{equation}\label{app_h:eval_a}
  E_i= -\frac{1}{\pi} \int_{\PT_i} q_{i}(y_1,y_2;\kappa_i) \,
  {\mathcal F}_{i}(y_1,y_2;\kappa_i)\, dy_1 dy_2  \,,
\end{equation}
where, in terms of a double integral over the patch, we have
\begin{equation}
  {\mathcal F}_{i}(y_1,y_2;\kappa_i)=\frac{1}{4\pi}\int_{\PT_i}
  q_{i}(y_1^{\p},y_{2}^{\p};\kappa_i)
  \log\left( \left(y_1-y_1^{\p}\right)^2+\left(y_2-y_{2}^{\p}\right)^2\right)\,
  dy_1^{\p} dy_2^{\p}\,.
\end{equation}
\esub
\end{lemma}

The proof of this result is analogous to that done in Appendix C of
\cite{GrebWard25} for the narrow capture MFPT problem {\em interior}
to a sphere, and is omitted.

As shown in Appendix C of \cite{GrebWard25}, when $\PT_i$ is a disk of
radius $a_i$ we can determine $E_i$ in (\ref{app_h:eval}) up to a
quadrature. For a locally circular patch, both the charge
density $q_i$ and the solution to (\ref{mfpt:fprob}) are radially symmetric in
$\rho_0=(y_1^2+y_2^2)^{1/2}$, and we obtain 
\begin{equation}\label{app_h:e_rad}
  E_i = - \int_{0}^{a_i} 2\rho_0 q_{i}(\rho_0;\kappa_i)
  {\mathcal F}_{i}(\rho_0;\kappa_i)
  \, d\rho_0 \,; \quad
 {\mathcal F}_{i,\rho_0} = \frac{1}{\rho_0} \int_{0}^{\rho_0} \eta
  q_{i}(\eta;\kappa_i) \, d\eta \,, \,\,\, 0\leq \rho_0\leq a_i\,,
\end{equation}
with ${\mathcal F}_{i} =\left({C_i/2}\right)\log{a_i}$ at $\rho=a_i$.
By integrating this result by parts we obtain (\ref{mfpt:Ej_all}).

To determine the limiting asymptotics in (\ref{mfpt:Ej_asy}) of Lemma
\ref{lemma:Ej_kappa} when $\PT_i$ is a disk, we proceed as in Appendix
C of \cite{GrebWard25}. When $\kappa_i=\infty$, we use
(\ref{mfpt:wc_q}) for $q_{i}(\rho_0;\infty)$ together with
$C_i=C_i(\infty)={2a_i/\pi}$ in (\ref{mfpt:Ej_all}). Upon evaluating
the resulting integrals analytically we obtain the expression for
$E_i(\infty)$ given in (\ref{mfpt:Ej_asy}) of Lemma
\ref{lemma:Ej_kappa}. Finally, to approximate $E_i$ for
$\kappa_i\ll 1$, we observe from (\ref{mfpt:wc}) that in this limit
$-\partial_{y_3}w_i\sim \kappa$ on the patch $y_3=0$,
$(y_1,y_2)\in \PT_i$. As a result, we identify that
$q_i(\rho_0;\kappa_i)\sim {\kappa_i/2}$ for $0\leq \rho_0\leq
a_i$. With this approximation for $q_i$, we can evaluate the integrals
in (\ref{mfpt:Ej_all}), while using $C_i\sim {\kappa_i a_i^2/2}$ for
$\kappa_i\ll 1$ from (\ref{mfpt:cj_small_b}). In this way, we readily
derive the limiting asymptotics for $E_i$ for $\kappa_i\ll 1$ as given
in (\ref{mfpt:Ej_asy}) of Lemma \ref{lemma:Ej_kappa}.

\vspace*{-0.2cm}
\bibliographystyle{plain}
\bibliography{references}

\end{document}